\documentclass[10pt]{amsart} 

\usepackage{amsmath}
\usepackage{amssymb}

\usepackage{graphicx}
\usepackage{xcolor}
\usepackage{subcaption}
\usepackage{multirow}
\usepackage{enumitem}
\usepackage{mathrsfs}
\usepackage[rightcaption]{sidecap}
\usepackage{rotating}
\usepackage{tikz}

\theoremstyle{definition}

\theoremstyle{remark}

\newcommand{\IC}{\text{IC}}
\newcommand{\BC}{\text{BC}}
\newcommand{\PDE}{\text{PDE}}
\newcommand{\SIS}{\text{SIS}}

\newcommand{\eps}{\varepsilon}

\title[A scaled TW-PINN for the KdV equation]{A scaled TW-PINN framework for soliton solutions of the Korteweg--de Vries equation}

\author{Seungwan Han}
\address{Department of Mathematics $\&$ POSTECH MINDS (Mathematical Institute for Data Science), Pohang University of Science and Technology, Pohang 37673, Korea}
\email{han97@postech.ac.kr}

\author{Daniel Olmos-Liceaga}
\address{Departamento de Matem\'aticas, Universidad de Sonora, Hermosillo, Sonora 83000, Mexico}
\email{daniel.olmos@unison.mx}

\author{Jiaxi Gu}
\address{Department of Mathematics $\&$ POSTECH MINDS (Mathematical Institute for Data Science), Pohang University of Science and Technology, Pohang 37673, Korea}
\email{jiaxigu@postech.ac.kr}

\author{Jae-Hun Jung}
\address{Department of Mathematics $\&$ POSTECH MINDS (Mathematical Institute for Data Science), Pohang University of Science and Technology, Pohang 37673, Korea}
\email{jung153@postech.ac.kr}

\makeatletter
\@namedef{subjclassname@2020}{\textup{}2020 Mathematics Subject Classification}
\makeatother

\begin{document}

\maketitle

\begin{abstract}
Accurate approximation of soliton solutions of the Korteweg--de Vries (KdV) equation remains challenging for physics-informed neural networks (PINNs) due to localized wave structures, high-order spatial derivatives, nonlinear interactions between solitons, and difficulties in long-time prediction.
In this work, we address these challenges using a scaling PINN framework to approximate KdV soliton solutions.
The framework transforms the governing equation into a normalized form, approximates its solution using a PINN solver, and reconstructs the solution in the original variables through inverse scaling.
We introduce scaling transformations and PINN solvers tailored to one- and two-soliton problems.
For the one-soliton case, we employ a PINN solver that incorporates a traveling-wave coordinate to capture soliton propagation.
For the two-soliton case, we employ a two-step training procedure for the PINN solver.
The first step learns the individual soliton profiles and estimates their propagation speeds, while the second step uses the estimated speeds to construct traveling-wave coordinates in the PINN solver and approximate the full two-soliton solution, including nonlinear interactions.
We also introduce soliton-informed sampling (SIS), which places additional collocation points around the soliton trajectories.
Numerical results show that the scaling PINN framework accurately captures one-soliton propagation, resolves two-soliton interactions, and preserves soliton profiles during long-time extrapolation in both cases.
An additional three-soliton experiment suggests the applicability of the proposed two-step training procedure to multi-soliton solutions.
\end{abstract}

\section{Introduction}
The Korteweg--de Vries (KdV) equation \cite{KdV1895,Miura1976,Drazin1989} is a canonical nonlinear dispersive equation arising in the study of long, weakly nonlinear waves in shallow water.
We consider the KdV equation of the form
\begin{equation}
\label{eq:KdV}
    u_t + \left(u^2\right)_x + \eps^2 u_{xxx} = 0,
\end{equation}
where $\eps$ is the dispersion parameter.
One fundamental analytical property of the KdV equation is that it admits a localized traveling wave solution whose profile is preserved during propagation.
The solitary wave, known as a soliton, has an amplitude, width, and propagation speed that are intrinsically related.
The integrable structure of the KdV equation further underlies the existence of multi-soliton solutions, in which individual solitons undergo nonlinear interactions and subsequently recover their original characteristics after separation, apart from phase shifts induced by the interactions \cite{Zabusky1965,Hirota1971,Drazin1989}.
This remarkable behavior has been studied by classical analytical techniques, including the inverse scattering transform \cite{Gardner1967}, the Lax pair formulation \cite{Lax1968}, and Hirota's bilinear method \cite{Hirota1971}.

Recently, machine learning methods have emerged as a promising approach for solving the KdV equation.
Various machine learning methods have been proposed, such as the symbolic calculation method \cite{Xie2025}, the hybrid machine learning methods \cite{Meuris2023,Pan2024}, operator learning methods \cite{Michalowska2025,Toscano2026,Nayak2026,Yu2026,Khodakarami2026}, and the physics-informed neural networks (PINNs) \cite{Raissi2019,Fang2022,Lin2023,Zhang2023,Song2024,Pu2024,Chen2024,Liu2024,Urban2025,Xiong2026,Wang2026}.
Specifically, we are interested in the numerical solutions of KdV solitons.
The constraints of conservation laws are incorporated into the loss function to improve the preservation of physical invariants \cite{Fang2022,Wang2026}.
Zhang et al. \cite{Zhang2023} introduced a symmetry-enhanced PINN where the invariant surface conditions induced by the Lie symmetries or non-classical symmetries of PDEs are embedded into the loss function.
In \cite{Pu2024}, Pu and Chen proposed the Lax pairs informed neural networks tailored for the integrable systems with Lax pairs.
Chen et al. \cite{Chen2024} investigated the influence of activation function selection on one- and two-soliton solutions.
Despite these advances, two issues remain relevant to parameterized KdV soliton problems.
First, when the dispersion parameter changes, PINNs generally need to be retrained.
Second, PINNs can have difficulty in maintaining predictive accuracy beyond the training time domain during long-time extrapolation \cite{Wang2025}.

To address these issues, we adopt the scaling PINN framework introduced in \cite{Scaled-PINN}.
In this framework, a scaling transformation converts the governing equation into a normalized form.
The solution of the normalized equation is then approximated using a PINN solver.
The approximate solution is transformed back to the original variables through inverse scaling, allowing solutions for different dispersion parameter values to be recovered without retraining.
Furthermore, the framework allows the PINN solver to be designed according to the characteristics of the solution.
We exploit this flexibility to design PINN solvers that improve long-time prediction by incorporating the propagation characteristics of solitons.
For the one-soliton problem, we employ a PINN solver that incorporates a traveling-wave coordinate to represent the propagation of a single soliton.
The main contribution of this study is the extension of this solver to interacting multi-soliton solutions, which cannot generally be represented by a single traveling-wave coordinate.
For the two-soliton problem, we introduce a two-step training procedure to capture the propagation and nonlinear interactions of the solitons.
Specifically, the first step learns the profiles and propagation speeds of individual solitons, while the second step uses the estimated speeds to construct traveling-wave coordinates and approximate the full two-soliton solution.
We also propose soliton-informed sampling (SIS), which generates additional collocation points near the soliton trajectories.

The remainder of this paper is organized as follows.
Section \ref{sec:Scaled-PINN} briefly reviews the scaling PINN framework introduced in \cite{Scaled-PINN}.
Section \ref{sec:1sol-PINN} presents its application to the one-soliton problem, together with numerical results.
Section \ref{sec:2sol-PINN} addresses the two-soliton problem using a two-step training procedure and the proposed SIS, followed by numerical results.
Section \ref{sec:3sol-PINN} examines the potential applicability of the two-step training procedure to the three-soliton problem.
Finally, Section \ref{sec:con} concludes the paper and discusses future research directions.

\section{Scaling PINN framework}
\label{sec:Scaled-PINN}
In this section, we briefly introduce the scaling PINN framework proposed in \cite{Scaled-PINN}, which was developed for traveling-wave solutions of the reaction--diffusion equation.
The framework first transforms the governing equation into a normalized form and then uses a PINN solver to approximate its solution.
The solution in the original variables is subsequently recovered through inverse scaling.

\begin{figure}
\centering
\includegraphics[width=0.8\columnwidth]{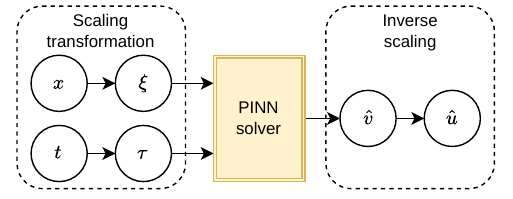}
\caption{Workflow of the scaling PINN framework, consisting of scaling transformation, PINN solver, and inverse scaling.}
\label{fig:Scaled-PINN}
\end{figure}

The overall framework is illustrated in Fig. \ref{fig:Scaled-PINN}, which summarizes three stages: (i) scaling transformation, (ii) solving the normalized equation using a PINN solver, and (iii) inverse scaling.
In \cite{Scaled-PINN}, a scaling transformation is applied to the reaction--diffusion equation to normalize the reaction and diffusion coefficients to unity.
This transformation maps the spatial and temporal coordinates $(x,t)$ to $(\xi,\tau)$.
The resulting normalized solution has a less steep profile, which makes the problem more suitable for PINN training.
A PINN solver is then used to obtain an approximate solution $\hat{v}$ of the normalized equation.
The resulting approximation is transformed back through inverse scaling to obtain an approximate solution $\hat{u}$ of the original equation.

In the present work, we apply the same framework to the KdV equation.
Both the scaling transformation and the PINN solver can be adapted to the governing equation and the desired solution characteristics.
We therefore introduce scaling transformations and PINN solvers designed for one- and multi-soliton solutions in Sec. \ref{sec:1sol-PINN} and Sec. \ref{sec:2sol-PINN}, respectively.
The resulting normalized solutions $\hat{v}$ are then transformed back to $\hat{u}$ through inverse scaling.

\section{Scaled TW-PINN for one-soliton}
\label{sec:1sol-PINN}
In this section, we apply the scaling PINN framework to the one-soliton solution of the KdV equation.
We first introduce the exact one-soliton solution and construct a scaling transformation to obtain the normalized equation.
We then describe the PINN solver designed for the one-soliton problem.
Finally, we present numerical results to evaluate the performance of the proposed solver.

\subsection{One-soliton solution and scaling}
The KdV equation \eqref{eq:KdV} admits a one-soliton solution that propagates as a traveling wave with a constant speed:
\begin{equation*}
    u(x,t) = A\,\text{sech}^{2} \left[ \frac{1}{\eps}\sqrt{\frac{A}{6}} \left(x-\frac{2A}{3}t \right) \right].
\end{equation*}
Here, $A$ denotes the soliton amplitude, and the corresponding propagation speed is
$c=\frac{2}{3}A$.
The characteristic spatial scale of the soliton is determined by both $A$ and $\eps$.
These relationships motivate the scaling transformation introduced below.

For a one-soliton solution, the scaling transformation is applied by incorporating both the amplitude $A$ and the dispersion parameter $\eps$.
The scaling transformation is defined as
\begin{equation}
\label{eq:1d_scale}
    t = \frac{\eps}{A\sqrt{A}}\tau, \qquad x = \frac{\eps}{\sqrt{A}}\xi, \qquad u = Av.
\end{equation}
Under this transformation, \eqref{eq:KdV} is converted into the normalized form
\begin{equation}
\label{eq:1d_normalize}
    v_{\tau} + (v^2)_{\xi} + v_{\xi\xi\xi}=0.
\end{equation}
The one-soliton solution of \eqref{eq:1d_normalize} is given by
\begin{equation*}
    v(\xi,\tau) = \text{sech}^2 \left( \frac{1}{\sqrt{6}} \left( \xi-\frac{2}{3}\tau \right) \right).
\end{equation*}
The normalized equation \eqref{eq:1d_normalize} is then solved using a PINN solver, and solutions with arbitrary values of $A$ and $\eps$ are recovered in the original variables through inverse scaling.

\subsection{One-soliton PINN solver and loss function}
\label{sec:loss}
\begin{figure}
\centering
\includegraphics[width=0.6\columnwidth]{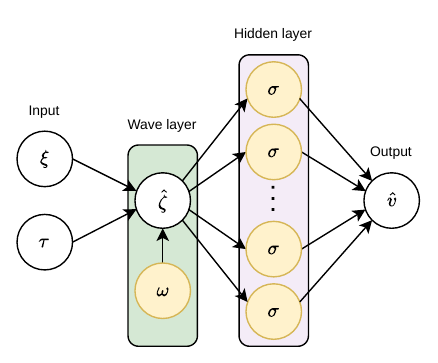}
\caption{Architecture of the TW-PINN for the one-soliton problem, with trainable components highlighted in yellow.}
\label{fig:model}
\end{figure}

The one-soliton solution presented above has a traveling-wave structure, preserving its shape during propagation.
This allows the normalized solution $v(\xi,\tau)$ to be expressed in terms of a single traveling-wave coordinate $\zeta$:
\begin{equation*}
    v(\xi,\tau)=V(\zeta), \qquad \zeta=\xi-c\tau,
\end{equation*}
where $c$ denotes the propagation speed of the traveling wave.
However, the standard PINN takes $(\xi,\tau)$ directly as inputs without explicit information about the soliton profile or propagation speed, and thus does not guarantee the traveling-wave structure.
To incorporate this structure into the network, we employ a wave layer that transforms $(\xi,\tau)$ into the traveling-wave coordinate.
In this layer, the propagation speed is represented by a trainable parameter $\omega$, and the traveling-wave coordinate is defined as
\begin{equation}
\label{eq:TW-coor}
    \hat{\zeta}=\xi-\omega\tau.
\end{equation}
The resulting coordinate $\hat{\zeta}$ is then passed to the hidden layers, which approximate the soliton profile.
We refer to the resulting network as a traveling-wave PINN (TW-PINN).
Fig. \ref{fig:model} illustrates the architecture of the TW-PINN used as the PINN solver for the one-soliton problem.

Based on this architecture, the TW-PINN approximates the solution $v(\xi,\tau)$ of \eqref{eq:1d_normalize} with the network prediction $\hat{v}(\xi,\tau)$.
The model is trained by minimizing a loss function that consists of the initial condition, boundary condition, and PDE residual losses.
The residual is defined as
\begin{equation}
\label{eq:resid}
    \mathcal{R}(\xi,\tau) := \hat{v}_{\tau} + (\hat{v}^{2})_{\xi} + \hat{v}_{\xi\xi\xi}.
\end{equation}
The total loss function is given by
\begin{equation}
\label{eq:loss}
    \mathcal{L} = \mathcal{L}_{\IC} + \mathcal{L}_{\BC} + \mathcal{L}_{\PDE},
\end{equation}
where the initial and boundary condition losses are defined as
\begin{align*}
    \mathcal{L}_{\IC} &= \frac{1}{N_{\IC}} \sum_{i=1}^{N_{\IC}} \left| \hat{v}(\xi_{\IC}^{i},\tau_{\IC}^{i}) - v_0(\xi_{\IC}^{i})\right|^2, \\
    \mathcal{L}_{\BC} &= \frac{1}{N_{\BC}} \sum_{i=1}^{N_{\BC}} \left| \hat{v}(\xi_{\BC}^{i},\tau_{\BC}^{i}) - v(\xi_{\BC}^{i},\tau_{\BC}^{i}) \right|^2,
\end{align*}
and the PDE residual loss is given by
\begin{equation*}
    \mathcal{L}_{\PDE} = \frac{1}{N_{\PDE}} \sum_{i=1}^{N_{\PDE}} \left| \mathcal{R} (\xi_{\PDE}^{i},\tau_{\PDE}^{i}) \right|^2.
\end{equation*}
Here, $\{v_0(\xi_{\IC}^i)\}_{i=1}^{N_{\IC}}$ denotes the prescribed initial condition, while $\{(\xi_{\BC}^i,\tau_{\BC}^i)\}_{i=1}^{N_{\BC}}$ denotes the prescribed boundary points.
The quantities $N_{\IC}$, $N_{\BC}$, and $N_{\PDE}$ denote the numbers of collocation points sampled from the initial condition, boundary condition, and the computational domain, respectively.

\subsection{Numerical results}
\label{sec:one_results}
For comparison, we consider the PINN with conservation laws proposed by Wang et al. \cite{Wang2026}, which was reported to preserve waveform structure during temporal evolution.
In this study, we refer to this method as PINN-CL.
The computational domain for both models is defined as $(\xi,\tau)\in[-15,20]\times[0,5]$.
The TW-PINN consists of a single hidden layer with $20$ neurons and uses the hyperbolic tangent activation function.
Latin hypercube sampling is employed to generate the collocation points.
The numbers of collocation points for the initial condition, boundary condition, and PDE residual are each set to $1{,}024$.
The TW-PINN is trained for $5{,}000$ iterations with the Adam optimizer, an initial learning rate of $0.01$, and a cosine annealing learning-rate schedule.
To verify the capability of the scaling to handle one-soliton solutions with different amplitudes and dispersion parameters, four representative parameter combinations are considered:
\begin{equation*}
    (A, \eps) = (1, 1),\, (1, 0.1),\, (4, 1),\, (4, 0.1).
\end{equation*}
Both models are trained once on the normalized equation and then used to recover one-soliton solutions for the different parameter combinations through the inverse scaling transformation.
We refer to the overall framework as the scaled PINN-CL or scaled TW-PINN, depending on whether the PINN with conservation laws or TW-PINN is used to solve the normalized equation.
During inference, the scaled coordinates $(\xi,\tau)$ corresponding to a given original coordinate $(x,t)$ are obtained from the scaling transformation in \eqref{eq:1d_scale}.
The recovered solutions are evaluated at $t=5$ and $10$.
The evaluation at $t=10$ corresponds to values of $\tau$ beyond the training domain $\tau\in[0,5]$, and larger $A$ or smaller $\eps$ leads to even larger values of $\tau$.
In this work, \textit{extrapolation} refers to cases in which the scaled inputs $(\xi,\tau)$ lie outside the training domain.

\begin{figure}
\centering
\includegraphics[width=0.85\columnwidth]{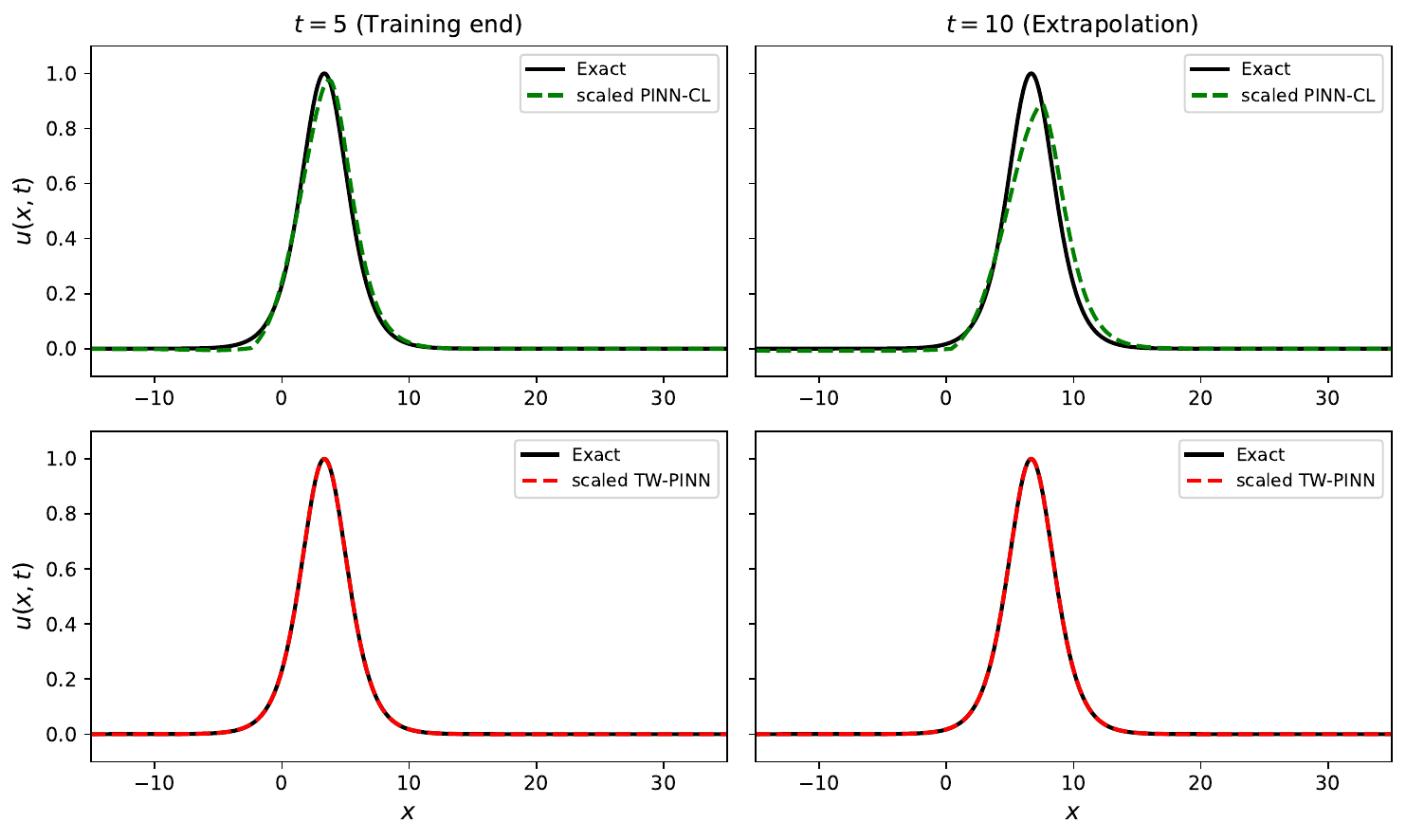}
\caption{Comparison of solution profiles obtained using the scaled PINN-CL (top) and the scaled TW-PINN (bottom) for the KdV equation with $A=1$ and $\eps=1$.}
\label{fig:eps1_a1}
\end{figure}

\begin{figure}
\centering
\includegraphics[width=0.85\columnwidth]{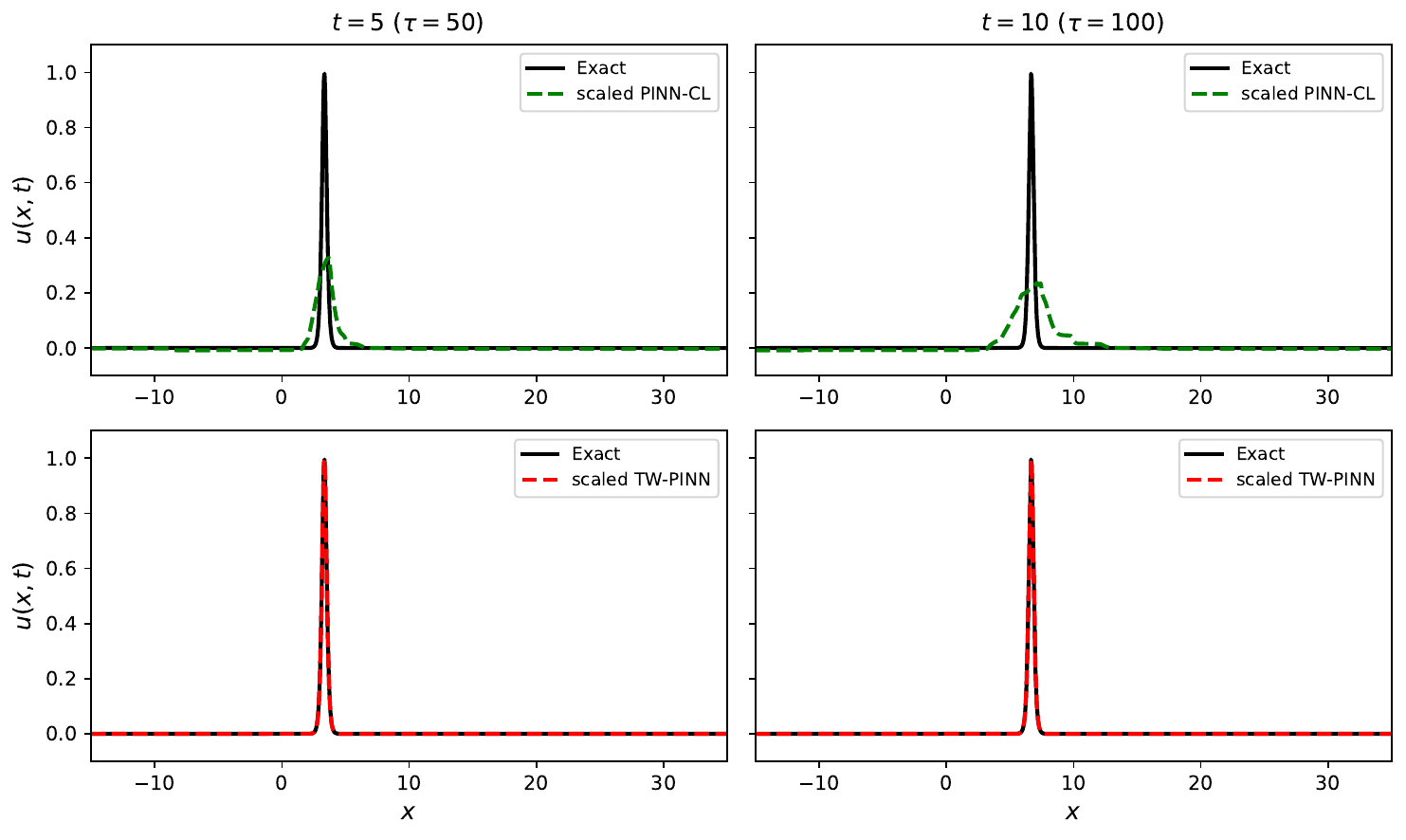}
\caption{Comparison of solution profiles obtained with the scaled PINN-CL (top) and scaled TW-PINN (bottom) for the KdV equation with $A=1$ and $\eps=0.1$.
The subplot titles show the original time $t$ and the corresponding scaled time $\tau$.}
\label{fig:eps0.1_a1}
\end{figure}

\begin{figure}
\centering
\includegraphics[width=0.85\columnwidth]{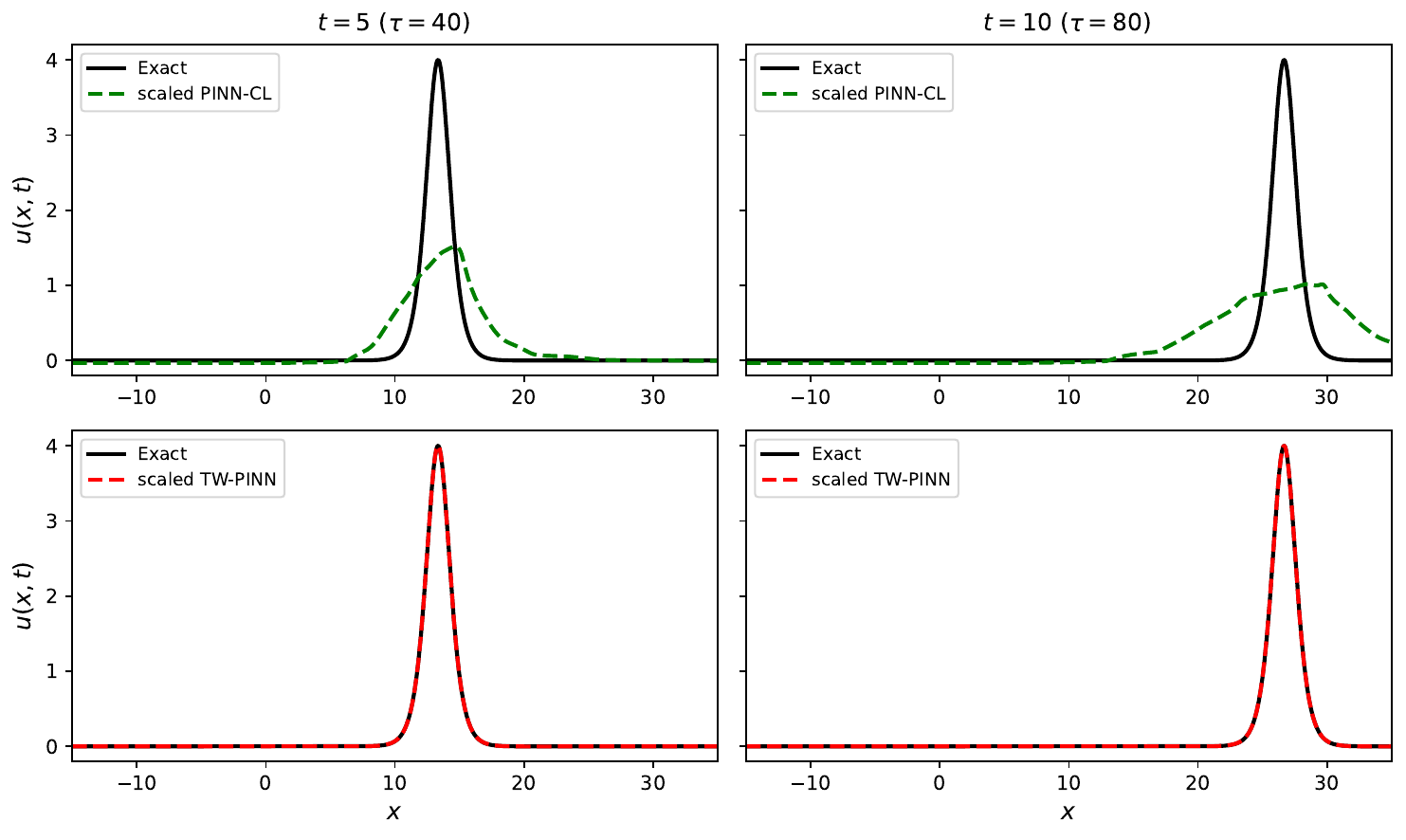}
\caption{Comparison of solution profiles obtained using the scaled PINN-CL (top) and the scaled TW-PINN (bottom) for the KdV equation with $A=4$ and $\eps=1$.
The subplot titles show the original time $t$ and the corresponding scaled time $\tau$.}
\label{fig:eps1_a5}
\end{figure}

\begin{figure}
\centering
\includegraphics[width=0.85\columnwidth]{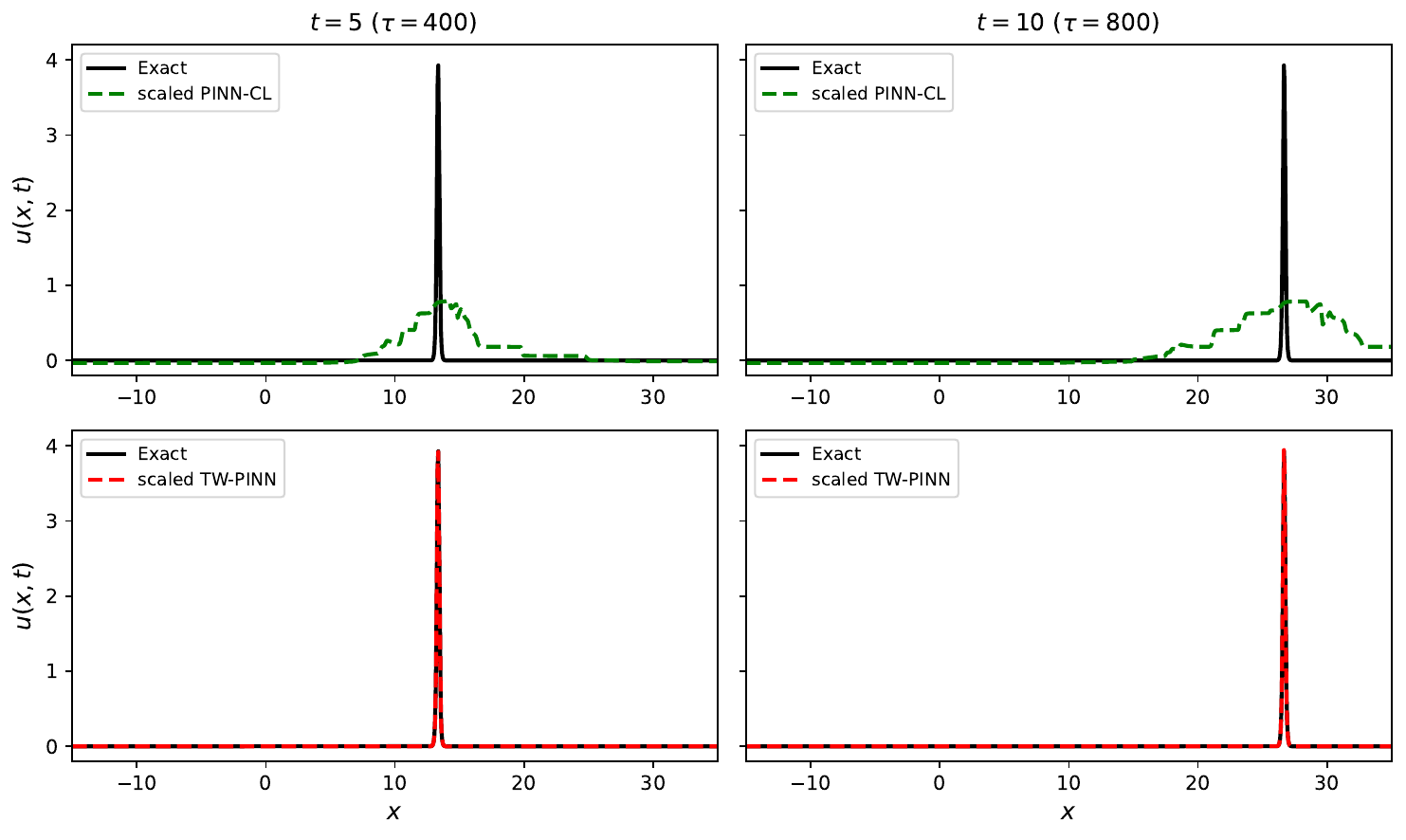}
\caption{Comparison of solution profiles obtained using the scaled PINN-CL (top) and the scaled TW-PINN (bottom) for the KdV equation with $A=4$ and $\eps=0.1$.
The subplot titles show the original time $t$ and the corresponding scaled time $\tau$.}
\label{fig:eps0.1_a5}
\end{figure}

Figures \ref{fig:eps1_a1}, \ref{fig:eps0.1_a1}, \ref{fig:eps1_a5}, and \ref{fig:eps0.1_a5} compare the exact one-soliton solutions with the predictions of the scaled PINN-CL and the scaled TW-PINN for the four representative parameter combinations.
All solution profiles are presented in the original coordinates $(x,t)$ after applying the inverse scaling transformation.

For $(A,\eps)=(1,1)$, both models successfully reproduce the overall one-soliton profile at $t=5$, for which the corresponding values of $\tau$ remain within the training domain.
At $t=10$, the scaled coordinates fall outside the training domain, where the scaled PINN-CL noticeably underestimates the peak amplitude, whereas the scaled TW-PINN continues to reproduce the soliton profile accurately.
For the remaining parameter combinations, $(1,0.1)$, $(4,1)$, and $(4,0.1)$, the performance difference between the two models becomes more pronounced, as the corresponding values of $\tau$ already exceed the upper bound of the training domain at $t=5$.
Beyond the training domain, the scaled PINN-CL exhibits increasingly pronounced oscillations and substantial underestimation of the peak amplitude, whereas the scaled TW-PINN maintains the soliton profile with high accuracy.
These results suggest that, for the parameter settings considered, the proposed TW-PINN provides more accurate long-time predictions than the scaled PINN-CL, which incorporates conservation laws into the loss function.
This advantage is particularly evident when extrapolation is required in the scaled coordinate space.

\section{Scaled TW-PINN for two-soliton}
\label{sec:2sol-PINN}
In this section, we extend the scaled TW-PINN framework to two-soliton solutions.
We first present the exact two-soliton solution and introduce a scaling transformation.
For two-soliton solutions, this transformation removes only the dispersion parameter because the individual soliton amplitudes cannot be normalized simultaneously.
In the PINN solver stage, the TW-PINN used for the one-soliton case cannot be directly applied because the two-soliton solution cannot be represented by a single traveling-wave coordinate.
We therefore introduce two traveling-wave coordinates and adopt a two-step training procedure.
In the first step, the propagation characteristics of the individual solitons are learned, while the second step incorporates this information into the two-soliton model through traveling-wave coordinates.
To improve training, we also introduce soliton-informed sampling, which generates additional collocation points around the soliton trajectories.
Finally, numerical experiments are presented to evaluate the performance of the proposed framework for the two-soliton problem.

The two-step training procedure is presented here for the two-soliton case, but it may also be applicable to general $N$-soliton solutions with appropriate modifications to the model architecture.
We examine this possibility through the three-soliton case in Sec. \ref{sec:3sol-PINN}, without modifying the training procedure.

\subsection{Two-soliton solution and scaling}
\label{sec:two-sol}
The two solitons generally propagate at different speeds, so the two-soliton solution cannot be represented by a single traveling-wave coordinate.
Instead, the two-soliton solution is expressed through the Hirota tau-function \cite{Hirota1971}.
Introducing the phase variables
\begin{equation*}
    \eta_1 = \frac{k_1x-k_1^3t+\delta_1}{\eps}, \qquad \eta_2 = \frac{k_2x-k_2^3t+\delta_2}{\eps},
\end{equation*}
where $k_i$ is the soliton parameter that determines the amplitude, width, and propagation speed of the $i$-th soliton, and $\delta_i$ is the phase constant.
The two-soliton solution is written as
\begin{equation*}
    u(x,t) = 6\eps^2 \frac{\partial^2}{\partial x^2} \left(\ln{F}\right),
\end{equation*}
where the tau-function is given by
\begin{equation*}
    F = 1+e^{\eta_1}+e^{\eta_2} + \alpha e^{\eta_1+\eta_2}.
\end{equation*}
Here, $\alpha$ is the interaction coefficient defined as
\begin{equation*}
    \alpha = \left(\frac{k_1-k_2}{k_1+k_2}\right)^2 .
\end{equation*}
The amplitude and propagation speed of each soliton are determined by its soliton parameter:
\begin{equation*}
    A_i = \frac{3}{2}k_i^2, \qquad c_i = k_i^2, \qquad i = 1,2.
\end{equation*}

The presence of two individual soliton parameters does not allow the scaling transformation to simultaneously normalize both amplitudes and the dispersion parameter.
We therefore introduce a scaling transformation that normalizes only the dispersion parameter $\eps$:
\begin{equation*}
    t = \eps \tau, \qquad x = \eps \xi, \qquad u=v.
\end{equation*}
Under this scaling transformation, \eqref{eq:KdV} becomes
\begin{equation*}
    v_{\tau} + (v^2)_{\xi} + v_{\xi\xi\xi}=0.
\end{equation*}
The phase variables are rewritten in the scaled coordinates as
\begin{equation*}
    \theta_1 = k_1 \xi - k_1^3 \tau + \frac{\delta_1}{\eps}, \qquad \theta_2 = k_2 \xi - k_2^3 \tau + \frac{\delta_2}{\eps}.
\end{equation*}
The corresponding two-soliton solution is
\begin{equation}
\label{eq:sol}
    v(\xi,\tau) = 6 \frac{\partial^2}{\partial \xi^2} \left(\ln F\right),
\end{equation}
where
\begin{equation*}
    F = 1 + e^{\theta_1}+e^{\theta_2} + \alpha e^{\theta_1+\theta_2}.
\end{equation*}

\begin{figure}
\centering
\includegraphics[width=\columnwidth]{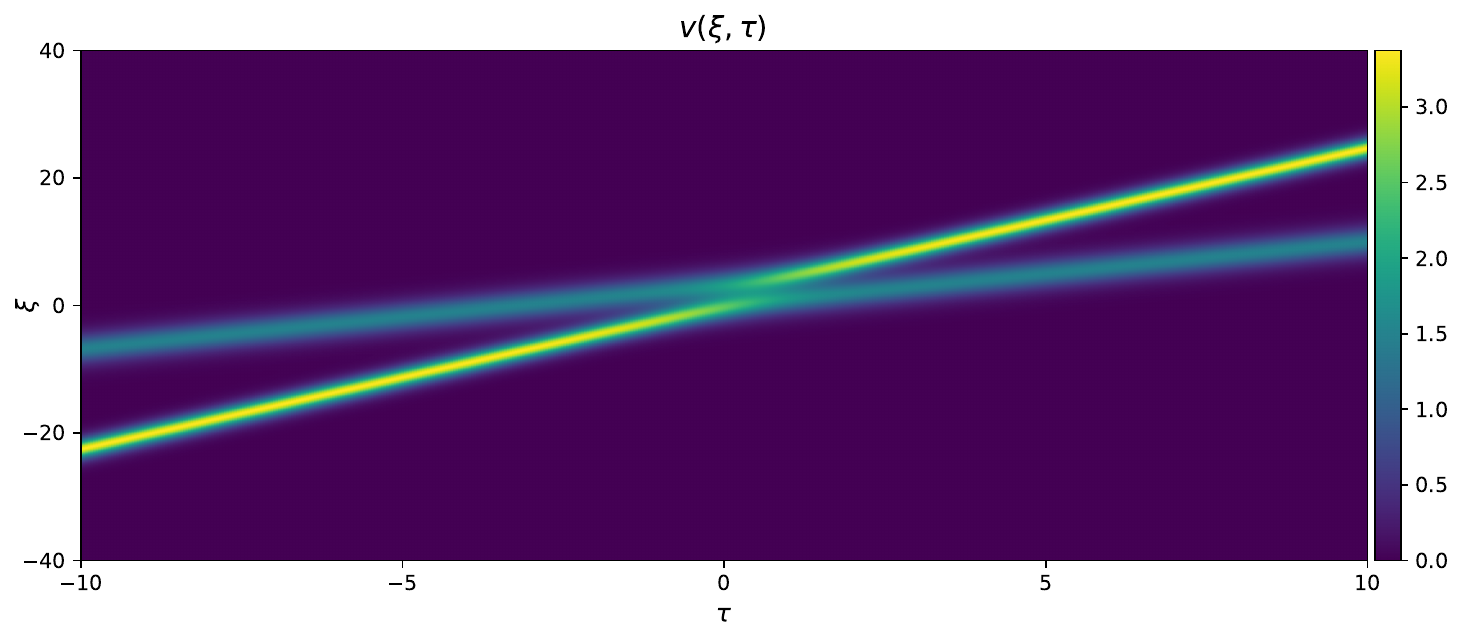}
\caption{Evolution of the exact two-soliton solution with $k_1=1.5$ and $k_2=1.0$, illustrating the propagation, nonlinear interaction, and subsequent separation of the two solitons.}
\label{fig:2-sol}
\end{figure}

The two-soliton solution exhibits nonlinear interactions between two localized waves with different profiles and propagation speeds.
Initially, the two solitons are well separated, and the faster soliton catches up with the slower one.
During the interaction, their profiles undergo temporary deformation.
After the interaction, the solitons separate and recover their original shapes and propagation speeds, with permanent phase shifts in their positions.
Fig. \ref{fig:2-sol} illustrates these dynamics for the soliton parameters $k_1=1.5$, $k_2=1.0$, $\delta_1=0$, and $\delta_2=0$ in \eqref{eq:sol}.

\subsection{Two-soliton PINN solver and soliton-informed sampling}
\begin{figure}
\centering
\includegraphics[width=\columnwidth]{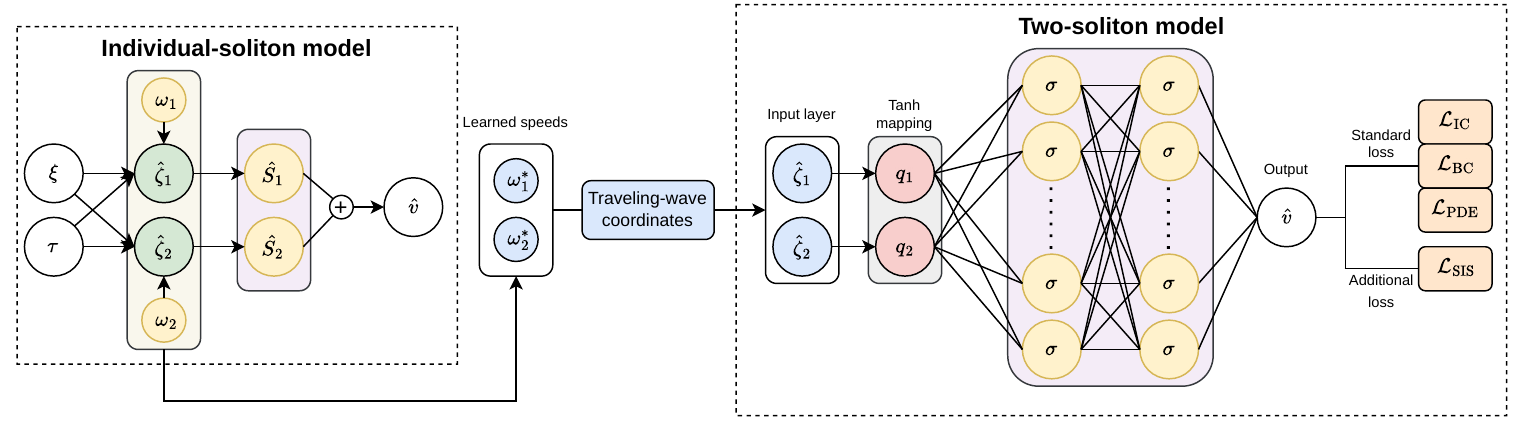}
\caption{Overview of the two-step training.
The first step learns the propagation characteristics of the individual solitons, and the second step incorporates this information through traveling-wave coordinates and uses the proposed SIS to improve the placement of collocation points during training for the full two-soliton solution.}
\label{fig:2steps}
\end{figure}

For the two-soliton problem, the TW-PINN is extended by introducing two traveling-wave coordinates to represent the propagation of the individual solitons.
Specifically, the network input is augmented with
\begin{equation*}
    \hat{\zeta}_1 = \xi - \omega_1\tau, \qquad \hat{\zeta}_2 = \xi - \omega_2\tau,
\end{equation*}
where $\omega_1$ and $\omega_2$ are intended to represent the propagation speeds of the two solitons.
Unlike in the one-soliton case, a single-step training procedure that directly optimizes $\omega_1$ and $\omega_2$ does not guarantee that the learned parameters correspond to the propagation speeds of the individual solitons.
For any distinct values of $\omega_1$ and $\omega_2$, the two traveling-wave coordinates $\hat{\zeta}_1$ and $\hat{\zeta}_2$ form an invertible linear transformation of the original inputs $(\xi,\tau)$.
The subsequent hidden layers can accommodate different choices of these parameters and still approximate the same two-soliton solution, so the actual propagation speeds may not be uniquely identified.
To address this difficulty in identifying the propagation speeds, we adopt a two-step training procedure in which the propagation characteristics of the individual solitons are first learned and then incorporated into the two-soliton model.
Fig. \ref{fig:2steps} illustrates the proposed two-step training procedure.

In the first step, the individual-soliton model is introduced to learn the propagation characteristics of the individual solitons before their nonlinear interaction.
Prior to the interaction, the two solitons are sufficiently separated, so that their mutual influence is negligible.
The solution can therefore be approximated by the superposition of two individual traveling-wave solutions:
\begin{equation*}
    v(\xi,\tau) \approx S_1(\xi-c_1\tau) + S_2(\xi-c_2\tau),
\end{equation*}
where $S_1$ and $S_2$ represent the individual soliton profiles, and $c_1$ and $c_2$ denote their corresponding propagation speeds.

\begin{figure}
\centering
\includegraphics[width=0.55\columnwidth]{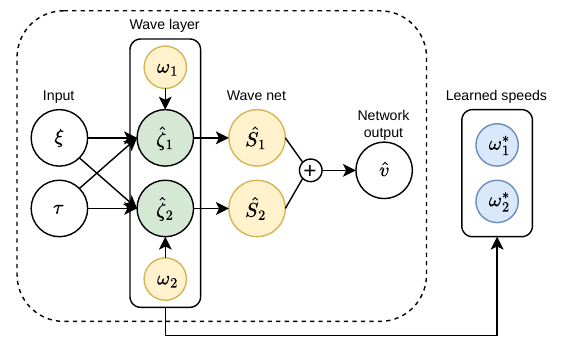}
\caption{Architecture of the individual-soliton model.
The trainable parameters $\omega_1$ and $\omega_2$ determine the propagation speeds of the individual solitons, while the subnetworks $\hat{S}_1$ and $\hat{S}_2$ learn their corresponding profiles.
After training, the optimized parameters $\omega_1^*$ and $\omega_2^*$ are extracted as the propagation speeds.}
\label{fig:step1}
\end{figure}

Based on this observation, the individual-soliton model is designed by exploiting the additive structure of the two-soliton solution.
Specifically, the model separately approximates the two individual soliton profiles, while the propagation speed of each soliton is treated as a trainable parameter.
Two traveling-wave coordinates are introduced as
\begin{equation*}
    \hat{\zeta}_i = \xi - \omega_i\tau, \qquad i = 1,2,
\end{equation*}
where $\omega_i$ denotes the trainable propagation speed of the $i$-th soliton.
The neural network approximation is then expressed as
\begin{equation*}
    \hat v(\xi,\tau) = \hat S_1(\hat{\zeta}_1) + \hat S_2(\hat{\zeta}_2),
\end{equation*}
where $\hat{S}_1$ and $\hat{S}_2$ denote the subnetworks representing the first and second soliton profiles, respectively.
Because each subnetwork is responsible for predicting only one individual soliton, the propagation parameter $\omega_i$ is learned together with the corresponding soliton profile.
The architecture of the individual-soliton model is illustrated in Fig. \ref{fig:step1}.
The inputs $(\xi,\tau)$ are transformed into the traveling-wave coordinates $\hat{\zeta}_1$ and $\hat{\zeta}_2$, which are then provided to the corresponding subnetworks.
The subnetworks have independent parameters and learn the two soliton profiles individually.
The two outputs are then summed to approximate the two-soliton solution when the solitons are well separated.

\begin{figure}
\centering
\includegraphics[width=0.8\columnwidth]{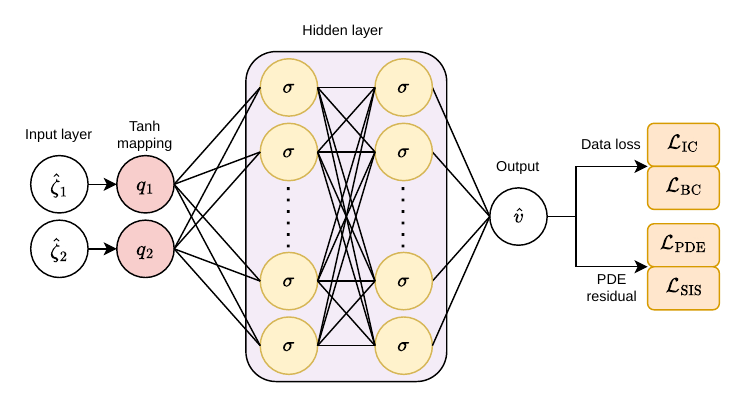}
\caption{Architecture of the two-soliton model.
The inputs are the traveling-wave coordinates $\hat{\zeta}_1$ and $\hat{\zeta}_2$, which are constructed using the propagation speeds learned by the individual-soliton model.}
\label{fig:step2}
\end{figure}

In the second step, the two-soliton model is trained to approximate the full solution, including the propagation of the individual solitons and their nonlinear interaction.
The propagation speeds learned by the individual-soliton model in Step 1 are used to define the traveling-wave coordinates for the two-soliton model as
\begin{equation*}
    \hat{\zeta}_i = \xi - \omega_i^*\tau, \qquad i = 1,2.
\end{equation*}
Here, $\omega_i^*$ denotes the propagation speed learned in Step 1 and is fixed throughout the training in Step 2.
Thus, $\hat{\zeta}_i$ varies with $(\xi,\tau)$, while the propagation speed $\omega_i^*$ remains fixed.
These coordinates are provided as inputs to the two-soliton model, allowing the network to incorporate the learned propagation characteristics of the individual solitons while learning their nonlinear interaction.
To facilitate long-time extrapolation, the traveling-wave coordinates are further mapped using the hyperbolic tangent function,
\begin{equation}
\label{eq:tanh}
    q_i = \tanh(\beta_i\hat{\zeta}_i), \qquad i = 1,2.
\end{equation}
Here, $\beta_i$ is a tunable factor that controls the steepness of the hyperbolic tangent mapping and its sensitivity to local variations in each traveling-wave coordinate.
We use the initial condition as a reference and apply constant shifts to the traveling-wave coordinates so that $q_i=0$ corresponds to the peak of each individual soliton at the initial time.
The hyperbolic tangent mapping reflects the localized nature of soliton profiles by emphasizing variations near the soliton centers, which may help reduce extrapolation errors.
The architecture of the two-soliton model is illustrated in Fig. \ref{fig:step2}.

In addition to the network architecture, the distribution of collocation points plays an important role in PINN training.
Conventional sampling methods distribute collocation points over the entire computational domain and may not provide sufficient resolution in regions where the solution varies rapidly.
To address this issue, adaptive sampling methods have been developed to concentrate more collocation points in regions that are difficult to resolve.
A representative approach is residual-based adaptive refinement \cite{Lu2021,Wu2023}, which selects additional collocation points by evaluating the PDE residual at candidate locations and refining regions with large residuals.

For soliton solutions, this issue becomes more pronounced because significant variations in the solution are concentrated around the moving soliton structures.
It is therefore important to allocate collocation points near these structures.
This issue has been addressed in recent PINN studies of the nonlinear Schr\"odinger equation, where adaptive sampling based on soliton proximity was employed to improve the accuracy of soliton solutions \cite{Moloshnikov2025}.
In this adaptive sampling strategy, the magnitude of the predicted solution is combined with the PDE residual to determine the sampling distribution.
Although this approach effectively exploits the localized structure of the solution, the sampling distribution depends on the current PINN prediction and is updated during training.
Thus, the placement of collocation points is governed by the predicted solution rather than by the underlying propagation characteristics of individual solitons.
Moreover, because the sampling distribution is determined from the magnitude of the current solution rather than from separately identified soliton trajectories, balanced sampling around individual solitons is not explicitly enforced in a multi-soliton solution.
This dependence may result in a suboptimal sampling distribution, particularly when the soliton structures have not yet been sufficiently resolved.

\begin{figure}
\centering
\includegraphics[width=0.9\columnwidth]{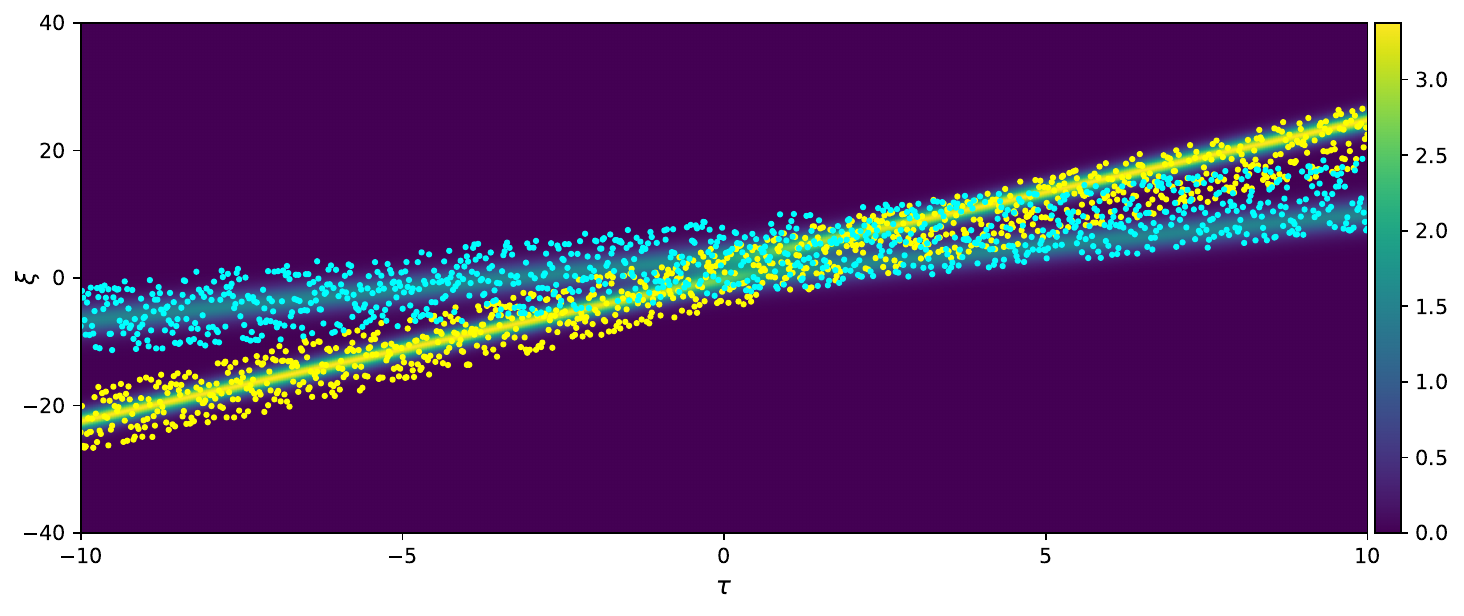}
\caption{Distribution of collocation points obtained using SIS over the exact two-soliton solution.
The yellow and cyan points denote the additional collocation points generated around the predicted trajectories of the first and second solitons, respectively.
The sampling radii are set to $r_1=5$ and $r_2=6$.}
\label{fig:sampling}
\end{figure}

To avoid this dependence, we propose \textit{soliton-informed sampling} (SIS), which uses soliton trajectories predicted from the speeds learned by the individual-soliton model to place collocation points for training the two-soliton model.
These trajectories are used to define sampling intervals around the individual solitons, within which additional collocation points are placed throughout the training time domain.
The sampling intervals are defined as follows.
Let $\xi_i(\tau)$ denote the predicted trajectory of the $i$-th soliton based on the speed learned by the individual-soliton model.
For each trajectory, a sampling interval
\begin{equation}
\label{eq:interval}
    [\xi_i (\tau) - r_i, \xi_i (\tau ) +r_i]
\end{equation}
is defined, where $r_i$ denotes the sampling radius associated with the $i$-th soliton.
The sampling radii $r_i$ are selected based on the soliton widths available from the initial condition, with an additional margin to account for possible positional displacements caused by phase shifts after the nonlinear interaction.
To generate the SIS collocation points, $N$ sampling times $\tau_1,\ldots,\tau_N$ are selected within the training time interval.
For each soliton, $N$ spatial coordinates are uniformly sampled from the corresponding sampling interval \eqref{eq:interval} at the selected sampling times.
Thus, $2N$ additional collocation points are generated for the two-soliton case.
Fig. \ref{fig:sampling} illustrates the resulting distribution of collocation points obtained using SIS with $N=1{,}024$.

The loss functions used to train the models in the two steps are defined as follows.
In Step 1, the individual-soliton model is trained using the same standard loss terms as the one-soliton model, namely the initial condition, boundary condition, and PDE residual losses defined in \eqref{eq:loss}.
In Step 2, the two-soliton model is trained using the same loss terms, together with an additional loss evaluated at the SIS collocation points:
\begin{equation*}
    \mathcal{L} = \mathcal{L}_{\IC} + \mathcal{L}_{\BC} + \mathcal{L}_{\PDE} + \mathcal{L}_{\SIS}.
\end{equation*}
The SIS loss is defined as
\begin{equation*}
    \mathcal{L}_{\SIS} = \frac{1}{N_{\SIS}}\sum_{i=1}^{N_{\SIS}}\left|\mathcal{R} (\xi_{\SIS}^{i},\tau_{\SIS}^{i})\right|^2,
\end{equation*}
where $\mathcal{R}$ denotes the PDE residual defined in \eqref{eq:resid}, and $N_{\SIS}$ is the total number of SIS collocation points.

\subsection{Numerical results}
\label{sec:two_results}
The numerical experiments evaluate the performance of the two-step training for the two-soliton problem.
The individual-soliton model in Step 1 is evaluated to verify whether the propagation characteristics of the individual solitons are accurately learned.
The two-soliton model in Step 2 is then evaluated to assess the effectiveness of the proposed SIS in improving prediction accuracy.
We compare three sampling strategies: baseline Latin hypercube sampling (LHS), LHS combined with residual-based adaptive refinement (RAR), and LHS combined with SIS.
Each model is trained ten times with different random seeds.
The figures show representative results from a single training run, whereas the quantitative results in the tables are reported as the mean and standard deviation over the ten independent runs.
The prediction errors are evaluated using the $L_2$ and $L_\infty$ errors, defined as
\begin{equation*}
 L_2 = \sqrt{\frac{1}{N} \sum_{i=1}^N \left( v^{(i)} - \hat{v}^{(i)} \right)^2}, \qquad
 L_{\infty} = \max_{1 \leqslant i \leqslant N} \left| v^{(i)} - \hat{v}^{(i)} \right|,
\end{equation*}
where $v^{(i)}$ denotes the exact solution and $\hat{v}^{(i)}$ denotes the predicted solution at the $i$-th of $N$ evaluation points.

The individual-soliton model is trained on $(\xi,\tau)\in[-35,10]\times[-11,-9]$, with the initial condition given at $\tau=-10$.
Since the purpose of Step 1 is to learn the propagation characteristics of the individual solitons before their nonlinear interaction, we select a short temporal interval around the initial time in which the two solitons remain sufficiently separated.

The standard PINN loss described in Sec. \ref{sec:loss} is used for training.
Each subnetwork consists of a single hidden layer with $20$ neurons and uses the hyperbolic tangent activation function.
The model is trained for $30{,}000$ iterations using the Adam optimizer with an initial learning rate of $0.01$ and a cosine annealing learning-rate schedule.
The numbers of collocation points for the initial condition, boundary condition, and PDE residual are each set to $1{,}024$, and the points are generated using Latin hypercube sampling.

\begin{figure}
\centering
\includegraphics[width=0.9\columnwidth]{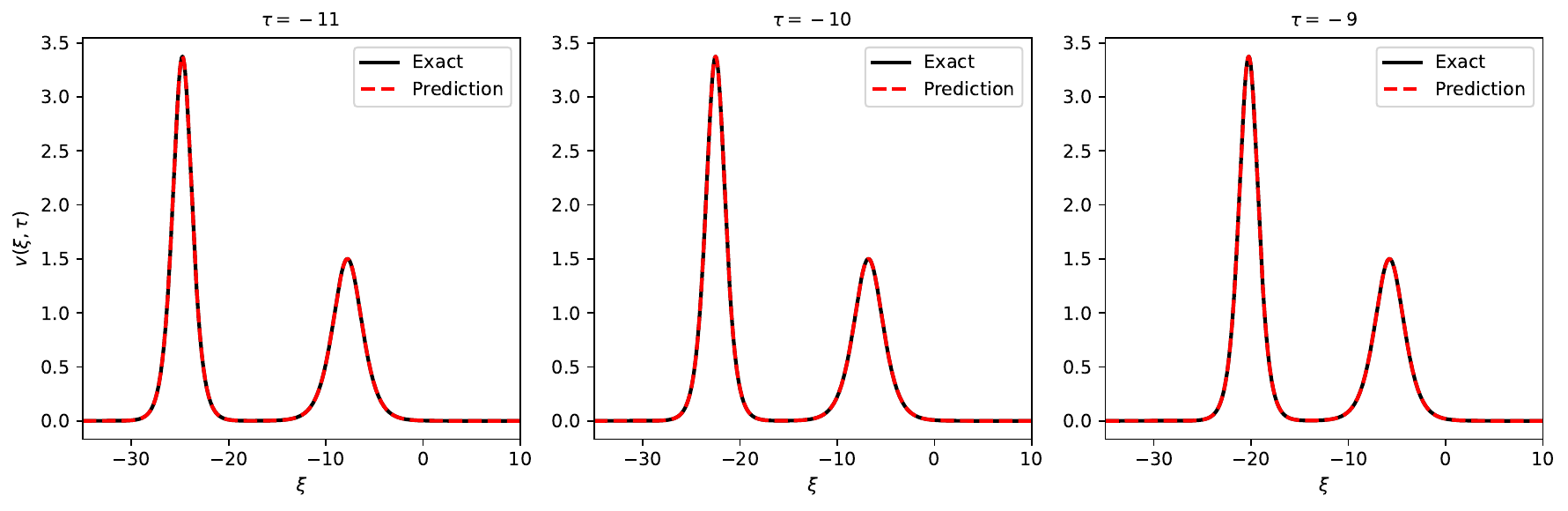}
\caption{Comparison of the exact and predicted soliton profiles at different time instances obtained using the individual-soliton model.}
\label{fig:M1}
\end{figure}

Fig. \ref{fig:M1} compares the predicted soliton profiles with the exact solutions at several times.
The predicted solutions closely match the exact profiles over the considered temporal interval, indicating that the individual soliton profiles are accurately reproduced by the individual-soliton model.

Table \ref{tab:M1} reports the mean absolute errors of the learned propagation speeds and the mean $L_2$ and $L_\infty$ errors of the predicted solutions, with the corresponding standard deviations given in parentheses.
The $L_2$ and $L_\infty$ errors are evaluated on a uniform $500\times500$ grid over the training domain.
The small errors in the learned propagation speeds and the predicted solutions indicate that the individual-soliton model accurately recovers the propagation characteristics required for Step 2.
\begin{table}
\centering
\caption{Absolute propagation speed errors with respect to the exact propagation speeds ($c_1=2.25$ and $c_2=1.0$), together with the $L_2$ and $L_{\infty}$ errors of the individual-soliton model.
Results are reported as the mean and standard deviation over $10$ independent runs.}
\label{tab:M1}
\resizebox{\textwidth}{!}{
\begin{tabular}{c c c c}
\hline
 $|c_1 - \omega_1|$  &  $|c_2 - \omega_2|$ & $L_2$ & $L_{\infty}$ \\ 
\hline
3.77e-03 (6.00e-03) & 2.83e-04 (5.13e-04) & 1.13e-03 (1.67e-03) & 9.16e-03 (1.38e-02) \\
\hline
\end{tabular}}
\end{table}

The two-soliton model is trained to predict the full two-soliton solution, including the propagation of the individual solitons and their nonlinear interaction.
The training domain in the $(\xi,\tau)$ coordinates is defined as $[-40,40]\times[-10,10]$.
In the numerical experiments, we use the exact propagation speeds to evaluate the use of traveling-wave coordinates and the sampling strategy without the influence of speed estimation errors from Step 1.
The parameters of the hyperbolic tangent mapping in \eqref{eq:tanh} are set to $\beta_1=\beta_2=0.2$.
These parameters are selected based on the spatial scales of the solitons observed from the initial condition and are fixed throughout training.
To investigate the effectiveness of the SIS, the model is trained using three sampling strategies, while the total number of PDE residual points is fixed at $3{,}072$.
The baseline uses $3{,}072$ globally distributed LHS points.
RAR uses $1{,}024$ LHS points and $2{,}048$ adaptively selected points from $10{,}000$ candidate points.
At every epoch, the PDE residual is evaluated at the candidate points, and the $2{,}048$ points with the largest residuals are selected.
The adaptive points from the previous epoch are replaced by the newly selected points.
SIS uses $1{,}024$ LHS points and two additional sets of $1{,}024$ PDE residual points.
Each set is distributed around a soliton trajectory determined by the corresponding propagation speed.
In all cases, $1{,}024$ points are used for the initial condition and $1{,}024$ points for the boundary condition.
The two-soliton model consists of two hidden layers with $120$ neurons per layer and uses the hyperbolic tangent activation function.
The model weights are initialized using Xavier initialization.
The model is first trained for $25{,}000$ iterations using the Adam optimizer with an initial learning rate of $0.01$ and a cosine annealing learning-rate schedule, followed by $5{,}000$ iterations using the L-BFGS optimizer.

\begin{figure}
\centering
\includegraphics[width=0.9\columnwidth]{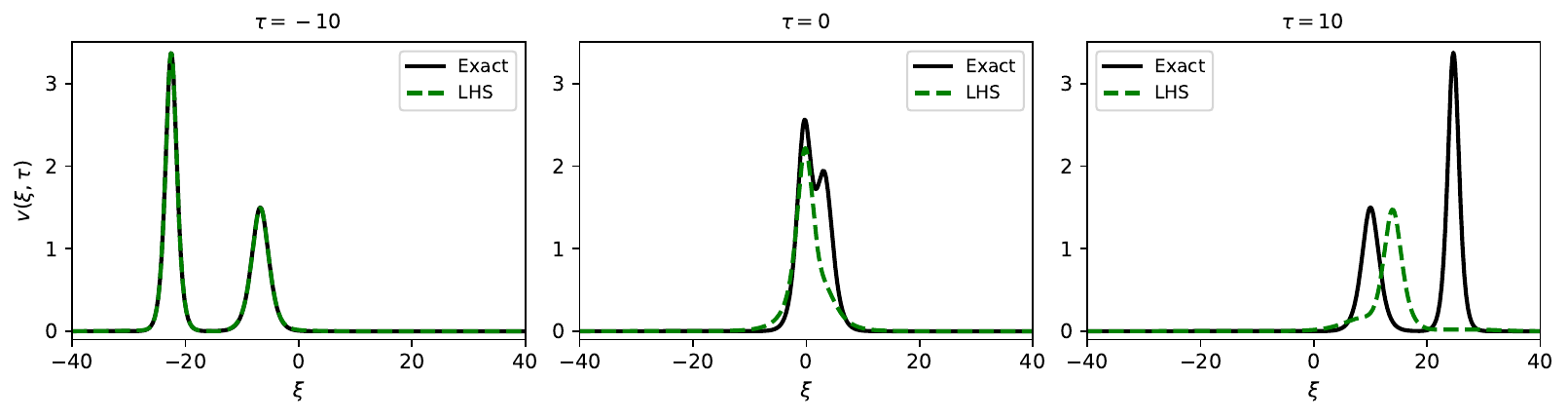}\\
\includegraphics[width=0.9\columnwidth]{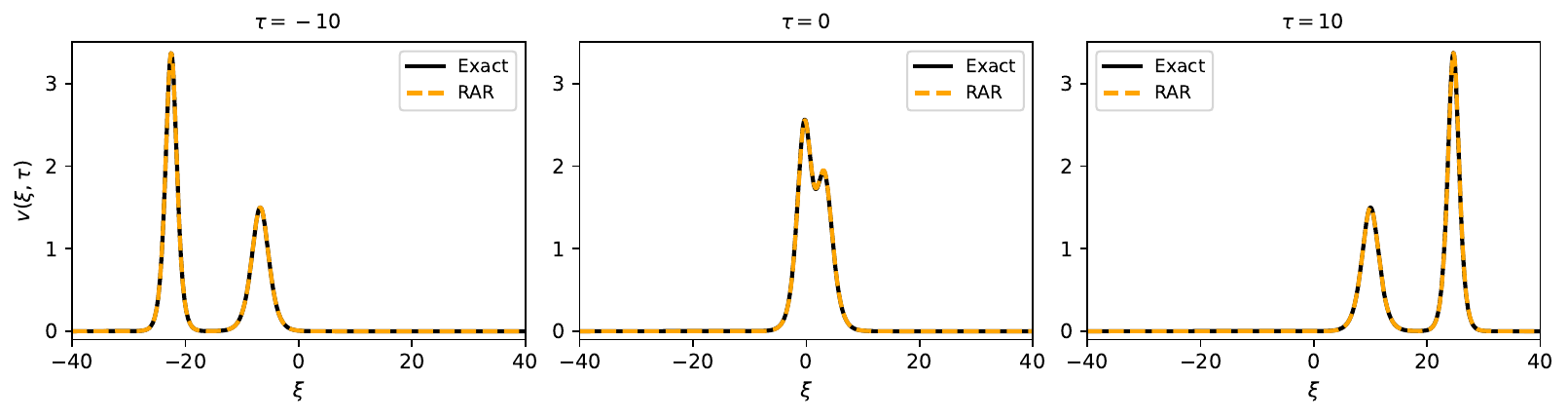}\\
\includegraphics[width=0.9\columnwidth]{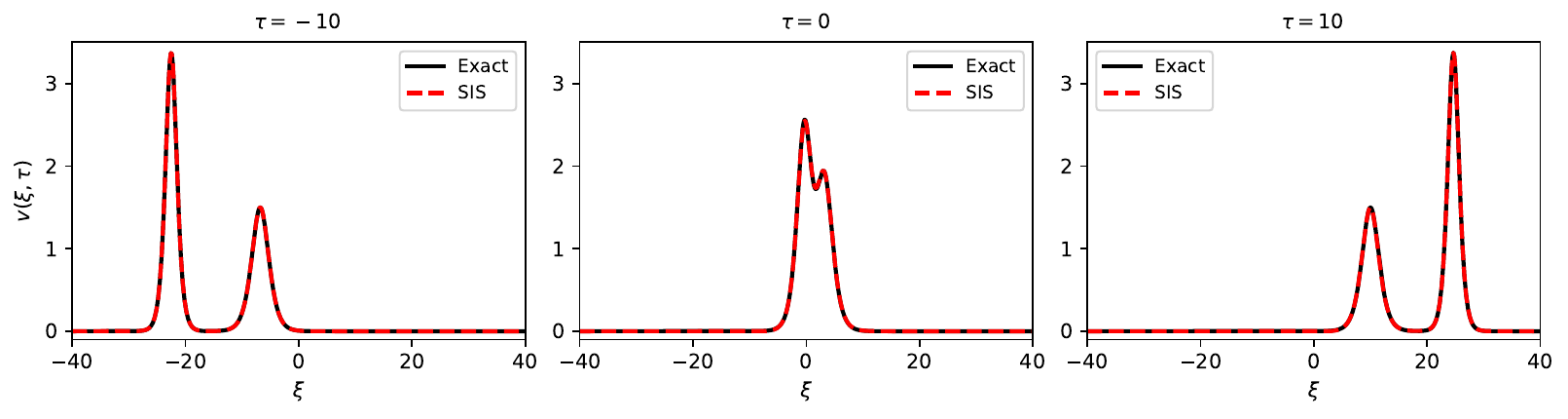}
\caption{Comparison of the exact and predicted two-soliton profiles obtained using LHS, RAR, and SIS at three representative time instances.
Rows correspond to the three sampling strategies, and columns represent the different time instances.}
\label{fig:M2}
\end{figure}

Fig. \ref{fig:M2} compares the predicted two-soliton profiles obtained using the three sampling strategies at three representative times, $\tau=-10$, $0$, and $10$, corresponding to the pre-interaction, interaction, and post-interaction phases, respectively.
For the model using LHS, the initial two-soliton profile is reasonably well reproduced at $\tau=-10$.
During the interaction at $\tau=0$, the two solitons merge into a single structure in the predicted solution and remain merged at $\tau=10$, resulting in a substantial deviation from the exact two-soliton solution.
In contrast, the predictions obtained using RAR and SIS closely match the exact solution throughout the evolution.

\begin{figure}
\centering
\includegraphics[width=0.9\columnwidth]{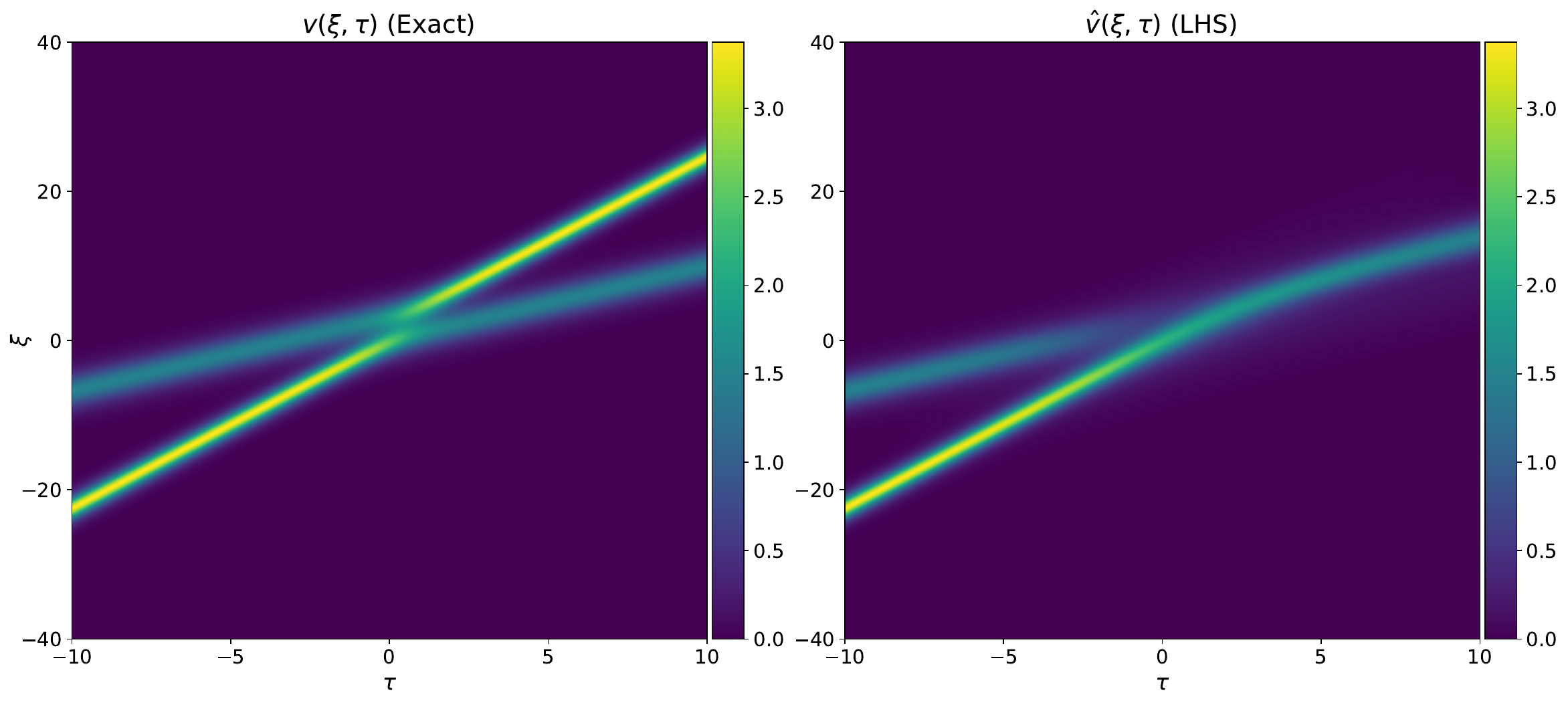}
\caption{Comparison of the exact two-soliton solution with the prediction of the model using LHS.}
\label{fig:LHS}
\end{figure}

To examine the LHS results over the entire training domain, the exact and predicted two-soliton solutions are compared through contour plots.
Fig. \ref{fig:LHS} shows the corresponding contours for the model using LHS.
The predicted solution initially captures the two individual solitons but fails to recover their separation after the nonlinear interaction, leading to a substantial deviation from the exact solution in the post-interaction region.

Since the model using LHS does not adequately preserve the two-soliton structure in our experiments, the following comparison focuses on RAR and SIS.
The quantitative results are summarized in Table \ref{tab:M2}.
The $L_2$ and $L_\infty$ errors are evaluated on a uniform $500\times500$ grid over the training domain.
To evaluate the phase shift $\Delta_i$ of each soliton, its peak positions at the initial and terminal times are identified from the corresponding local maxima on a uniform spatial grid with $d\xi=0.01$.
Each peak position is then refined by quadratic interpolation using the peak grid point and its two neighboring points.
The absolute errors in the estimated phase shifts are computed relative to the exact values.
All error metrics are reported as means and standard deviations over independent training runs, with the standard deviations given in parentheses.
SIS yields lower errors than RAR across all evaluation metrics, including the $L_2$ and $L_{\infty}$ errors and the phase-shift errors of both solitons.
In addition, RAR requires repeated residual evaluations over a large set of candidate points to identify regions for refinement, whereas SIS uses the soliton trajectories obtained from the individual-soliton model to place additional collocation points around the moving soliton structures.
These results indicate that SIS provides an effective strategy for collocation-point placement and offers an alternative to residual-based refinement.
\begin{table}
\centering
\caption{Quantitative comparison of prediction and phase-shift errors for the two-soliton model using RAR and SIS with the same total number of PDE residual points.
Results are reported as the mean and standard deviation.}
\label{tab:M2}
\resizebox{\textwidth}{!}{
\begin{tabular}{l c c c c}
\hline
Sampling strategy&   $L_2$  &  $L_{\infty}$ & $\left| \Delta_1^{\text{exact}} - \Delta_1^{\text{pred}}\right|$ & $\left| \Delta_2^{\text{exact}} - \Delta_2^{\text{pred}}\right|$ \\ 
\hline
RAR           & 8.03e-03 (4.20e-03) & 7.85e-02 (5.81e-02) & 3.66e-02 (4.21e-02) & 1.98e-02 (1.75e-02) \\
SIS           & 4.18e-03 (2.56e-03) & 4.23e-02 (2.78e-02) & 1.94e-02 (1.65e-02) & 1.08e-02 (1.15e-02) \\
\hline
\end{tabular}}
\end{table}

\begin{figure}
\centering
\includegraphics[width=0.9\columnwidth]{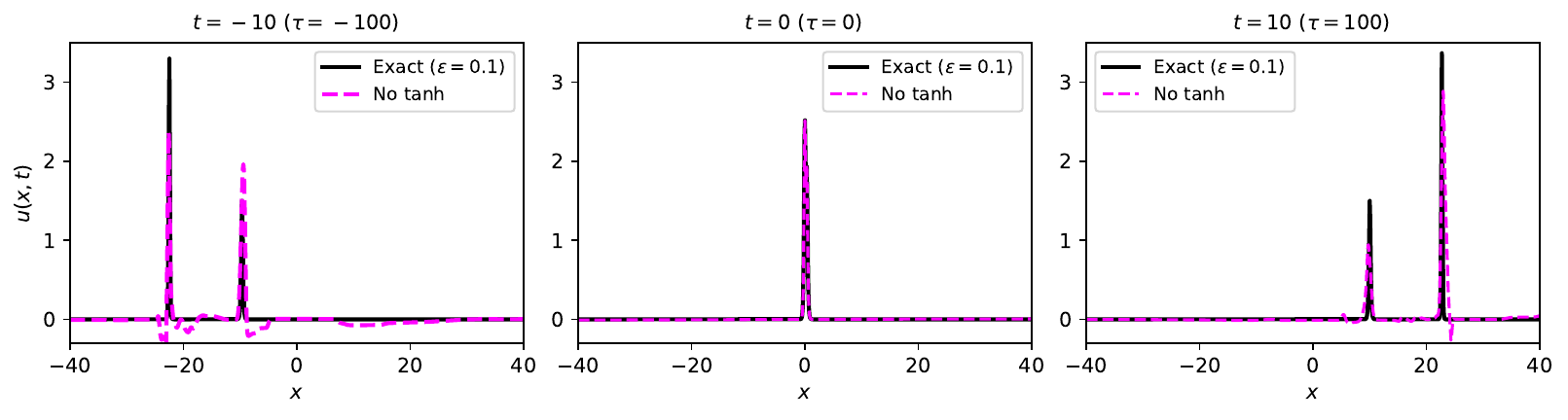}\\
\includegraphics[width=0.9\columnwidth]{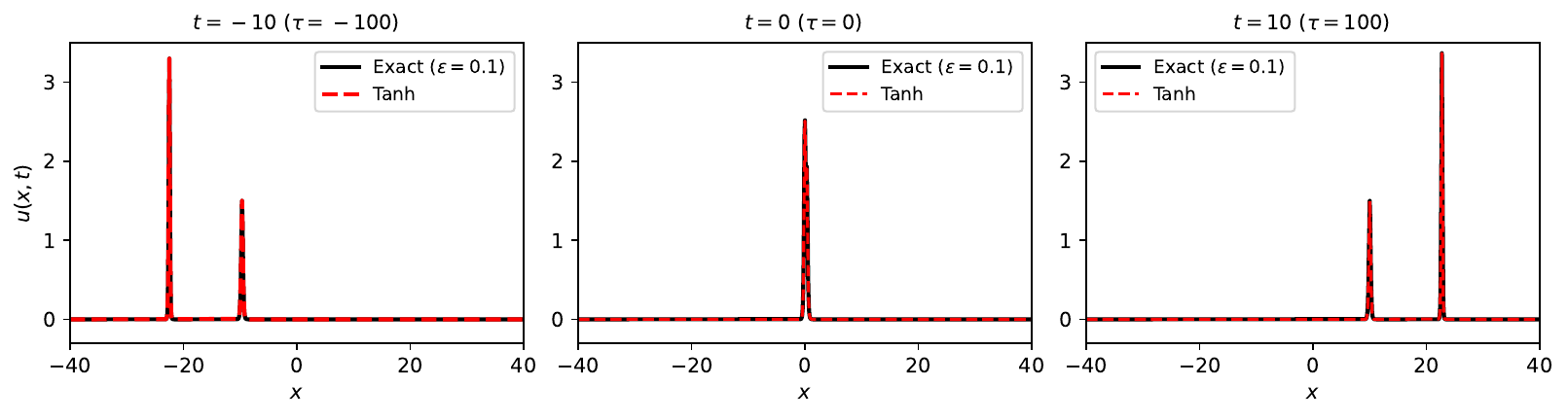}
\caption{Exact and predicted two-soliton solution profiles for $\eps=0.1$ at $t=-10$, $0$, and $10$. Results without and with the tanh mapping are compared to assess long-time extrapolation beyond the training domain.}
\label{fig:eps}
\end{figure}

\begin{figure}
\centering
\includegraphics[width=0.6\columnwidth]{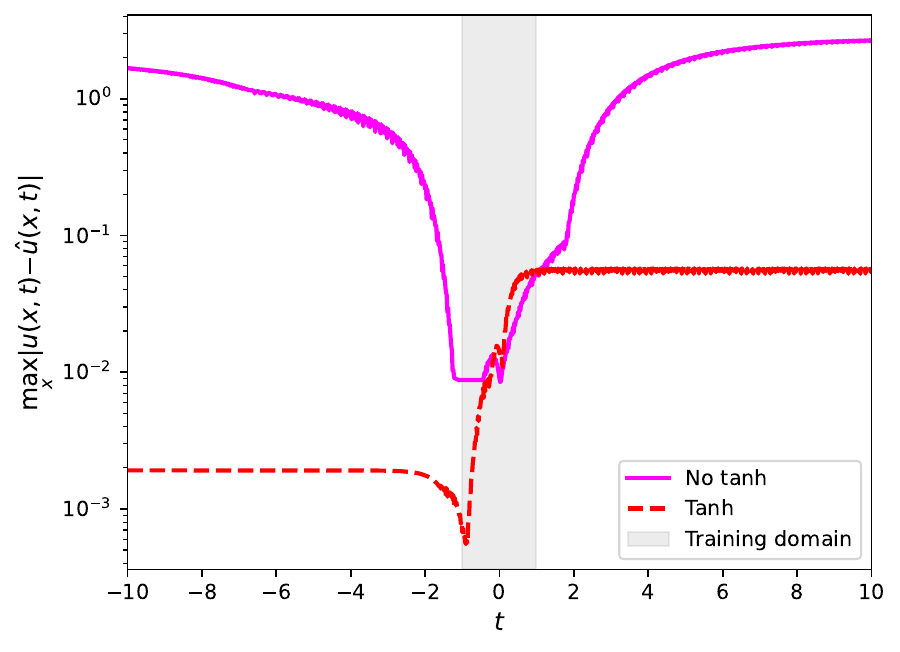}
\caption{Maximum spatial error between the exact and predicted two-soliton solutions over time for the two-soliton model without and with the tanh mapping.
The shaded region denotes the training domain.}
\label{fig:max}
\end{figure}

We next examine the long-time extrapolation capability of the two-soliton model beyond the training domain by comparing models with and without the tanh mapping.
Both models use SIS with identical model and training settings, differing only in whether the tanh mapping is applied.
This capability is particularly important for the scaling PINN framework because accurate reconstruction in the original variables requires reliable predictions over an extended range of scaled coordinates.
We therefore consider $\eps=0.1$, for which the scaling transformation leads to an expanded range of scaled coordinates.
The evaluation domain is set to $(\xi,\tau)\in[-400,400]\times[-100,100]$, which is ten times larger than the training domain $(\xi,\tau)\in[-40,40]\times[-10,10]$ along each coordinate.
Fig.~\ref{fig:eps} shows the exact and predicted two-soliton profiles at $t=-10$, $0$, and $10$, while Fig.~\ref{fig:max} presents the maximum spatial errors over the entire evaluation time interval.

The solution profiles in Fig. \ref{fig:eps} show a clear difference between the two models.
Without the tanh mapping, the model predicts the soliton positions reasonably well, but noticeable discrepancies appear in the soliton profiles and spurious oscillations develop outside the soliton regions.
In contrast, the model with the tanh mapping accurately reproduces the exact solution profiles even far beyond the training domain.

The difference is also evident in the maximum errors shown in Fig. \ref{fig:max}.
Without the tanh mapping, the maximum error increases progressively as the evaluation time moves farther beyond the training domain.
With the tanh mapping, the maximum error remains relatively small and does not increase substantially even outside the training domain.

These results demonstrate that the tanh mapping improves the long-time extrapolation capability of the two-soliton model.
This improvement allows accurate reconstruction of the solution in the original variables through inverse scaling, even when the scaled coordinates extend beyond the training domain.

\subsection{Conservation laws}
In Sec. \ref{sec:one_results}, we showed that TW-PINN, which incorporates the traveling-wave structure into the neural network architecture, better preserved the soliton structure during long-time extrapolation than PINN-CL for the one-soliton solution.
Since these two approaches act on different components of a PINN, namely the network architecture and the loss function, they can be applied simultaneously.
We therefore investigate whether the addition of conservation laws further improves the prediction accuracy of TW-PINN for the two-soliton solution.
To this end, we compare three configurations: (i) TW-PINN with SIS, (ii) TW-PINN with both SIS and conservation laws, and (iii) TW-PINN with conservation laws but without SIS.
The third configuration uses the same LHS setting as in Sec.~\ref{sec:two_results}, with conservation laws added to the loss function.
It serves as a reference for assessing whether conservation laws can achieve prediction accuracy comparable to that obtained with SIS.

Table \ref{tab:conservation} presents the $L_2$ and $L_\infty$ prediction errors and the absolute errors in mass and momentum conservation for the three configurations.
The addition of conservation laws to the SIS configuration slightly reduces the $L_2$ and $L_\infty$ prediction errors, whereas the errors in mass and momentum conservation decrease substantially.
The configuration with conservation laws but without SIS also exhibits relatively small conservation errors, but its $L_2$ and $L_\infty$ prediction errors are larger than those of the SIS configurations.
This difference may be attributed to the distinct roles of SIS and conservation laws. SIS introduces additional training points in the soliton regions and thereby places greater emphasis on these regions during training, whereas conservation laws impose constraints on global quantities without directly targeting local errors.
\begin{table}
\centering
\caption{Comparison of the $L_2$ and $L_\infty$ prediction errors and the absolute errors in mass and momentum for the three configurations: SIS, SIS with conservation laws, and conservation laws without SIS.
Values are reported as means, with standard deviations in parentheses.
Mass and momentum errors are calculated as absolute differences from their exact values of 15 and 26.25, respectively.
}
\label{tab:conservation}
\resizebox{\textwidth}{!}{
\begin{tabular}{l c c c c}
\hline
Additional loss terms &    $L_2$  & $L_{\infty}$ & Mass error&  Momentum error\\ 
\hline
SIS  & 4.18e-03 (2.56e-03) & 4.23e-02 (2.78e-02) &  1.00e-02 (3.91e-03) & 2.08e-02 (8.92E-03)\\
SIS + Conservation & 3.77e-03 (1.36e-03) & 4.07e-02 (1.64e-02) & 3.06e-03 (1.56e-03) & 6.09e-03 (3.02E-03)\\
Conservation  & 4.31e-02 (2.48e-02) & 4.11e-01 (2.32e-01) & 6.80e-03 (3.04e-03) & 8.31e-03 (4.45e-03)\\
\hline
\end{tabular}}
\end{table}

\section{Scaled TW-PINN for multi-soliton}
\label{sec:3sol-PINN}
In this section, we apply the two-step training procedure to the three-soliton problem to examine its applicability to multi-soliton solutions.
We first present the exact three-soliton solution in terms of the Hirota tau-function.
We then modify the model architecture by introducing an additional traveling-wave coordinate for the third soliton, without changing the training procedure.
Finally, we present numerical results to evaluate the performance of the resulting model.

\subsection{Three-soliton solution}
The Hirota tau-function provides an exact representation of $N$-soliton solutions of the KdV equation \cite{Hirota1971}.
The two-soliton solution considered in Sec. \ref{sec:two-sol} corresponds to the $N=2$ case of this formulation.
For the three-soliton solution, introducing the phase variables
\begin{equation*}
    \eta_i = \frac{k_i x-k_i^3 t + \delta_i}{\eps}, \qquad i = 1,2,3,
\end{equation*}
the exact solution is written as
\begin{equation*}
    u(x,t) = 6\eps^2 \frac{\partial^2}{\partial x^2}\left(\ln F\right),
\end{equation*}
where the tau-function is given by
\begin{equation*}
    F = 1+\sum_{i=1}^{3}e^{\eta_i} + \sum_{1\leq i<j\leq3}\alpha_{ij}e^{\eta_i+\eta_j} + \alpha_{12}\alpha_{13}\alpha_{23}e^{\eta_1+\eta_2+\eta_3},
\end{equation*}
with
\begin{equation*}
    \alpha_{ij} = \left(\frac{k_i-k_j}{k_i+k_j}\right)^2.
\end{equation*}
Here, $k_i$ determines the amplitude and propagation speed of the $i$-th soliton, while $\delta_i$ denotes its phase constant.
The pairwise interaction coefficients $\alpha_{ij}$ account for the interactions between each pair of solitons, while the final term represents the contribution of the three-soliton interaction.

For the three-soliton solution, we use the same scaling transformation as in Sec. \ref{sec:two-sol}.
The phase variables in the scaled coordinates are given by
\begin{equation*}
    \theta_i = k_i \xi - k_i^3 \tau + \frac{\delta_i}{\eps}, \qquad i=1,2,3.
\end{equation*}
The corresponding three-soliton solution is
\begin{equation*}
    v(\xi,\tau)
    = 6 \frac{\partial^2}{\partial \xi^2}
    \left(\ln F\right),
\end{equation*}
where
\begin{equation*}
    F = 1 + \sum_{i=1}^{3} e^{\theta_i}
    + \sum_{1 \le i < j \le 3}
    \alpha_{ij} e^{\theta_i+\theta_j}
    + \alpha_{12}\alpha_{13}\alpha_{23}
    e^{\theta_1+\theta_2+\theta_3}.
\end{equation*}

\subsection{Extension of scaled TW-PINN}
\begin{figure}
\centering
\includegraphics[width=0.9\columnwidth]{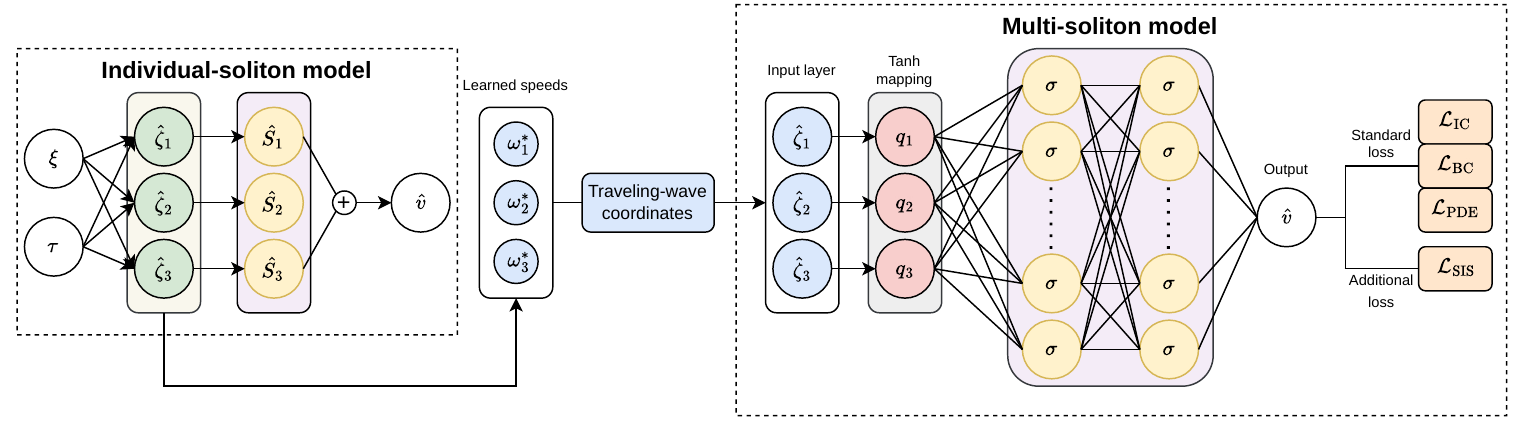}
\caption{Application of the two-step training procedure to the three-soliton problem through an extended model architecture.}
\label{fig:3steps}
\end{figure}

The two-step training procedure can be readily applied to the three-soliton problem without modification by extending the model architecture to include an additional traveling-wave coordinate for the third soliton.
The resulting architecture is illustrated in Fig. \ref{fig:3steps}.
In Step 1, the individual-soliton model employs the additional coordinate $\hat{\zeta}_3$ to learn the three individual soliton profiles, $\hat{S}_1$, $\hat{S}_2$, and $\hat{S}_3$, together with their propagation speeds.
The resulting information is then incorporated into Step 2 through the three traveling-wave coordinates $\hat{\zeta}_1$, $\hat{\zeta}_2$, and $\hat{\zeta}_3$.
The SIS strategy is similarly extended by placing three sets of additional collocation points around the predicted trajectories of the individual solitons.

\subsection{Numerical results}
We consider the exact three-soliton solution with $k_1=1.8$, $k_2=1.5$, and $k_3=1.0$, with $\delta_1=\delta_2=\delta_3=0$.
We use the same hidden-layer configuration as the two-soliton model, consisting of two hidden layers with $120$ neurons each.
The training domain is $(\xi,\tau)\in[-50,50]\times[-10,10]$.
For the initial condition, boundary condition, and PDE residual, we use $1{,}024$ collocation points each.
SIS generates three additional sets of $1{,}024$ collocation points, one around each soliton trajectory.
The sampling radii are $r_1=6$, $r_2=8$, and $r_3=10$, while the tanh mapping parameters are $\beta_1=\beta_2=\beta_3=0.15$.
The total training loss includes the initial condition, boundary condition, PDE residual, and SIS losses.

\begin{figure}
\centering
\includegraphics[width=0.9\columnwidth]{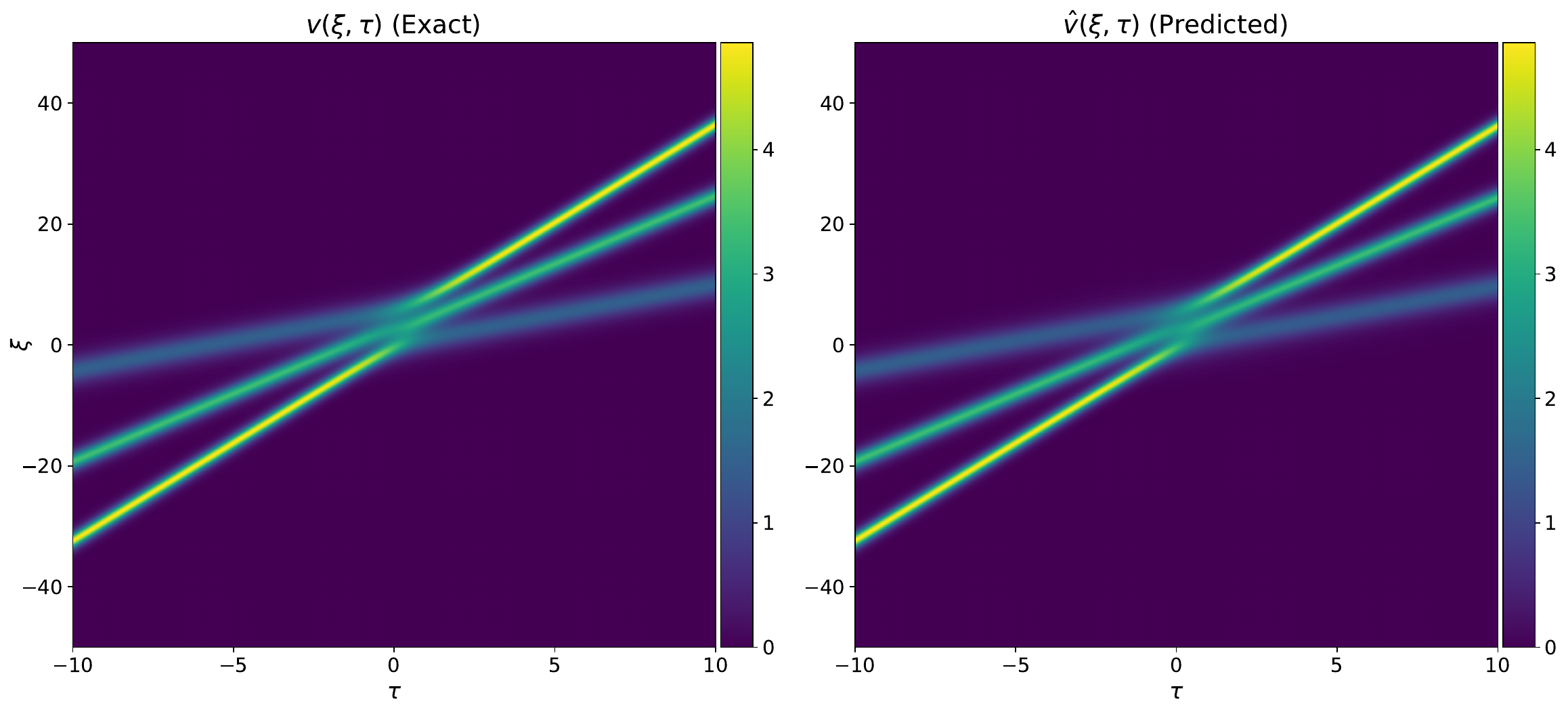}
\caption{Exact and predicted three-soliton solutions obtained using the multi-soliton model.}
\label{fig:3sol}
\end{figure}

\begin{figure}
\centering
\includegraphics[width=0.9\columnwidth]{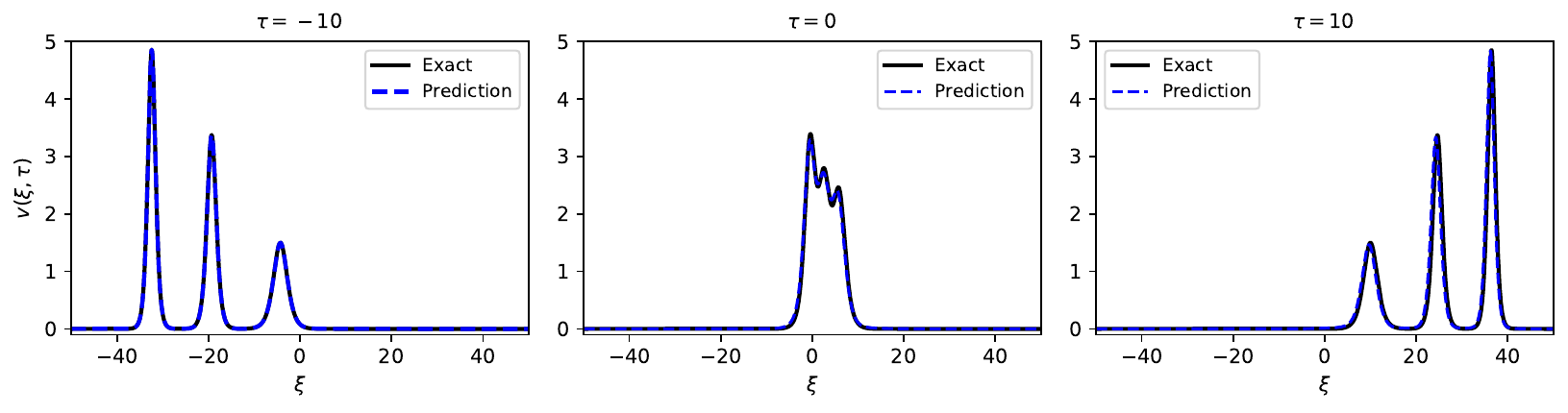}
\caption{Comparison of the exact and predicted three-soliton profiles at representative time instances.}
\label{fig:3sol_profile}
\end{figure}

Fig. \ref{fig:3sol} compares the exact and predicted three-soliton solutions.
The predicted solution preserves the three-soliton structure before and after the nonlinear interaction and accurately reproduces their interaction and subsequent separation.
The corresponding one-dimensional profiles at three representative times before, during, and after the nonlinear interaction are presented in Fig. \ref{fig:3sol_profile}.
The predicted profiles remain in close agreement with the exact profiles throughout the evolution, including during the nonlinear interaction and subsequent propagation.
These results show that the multi-soliton model captures the propagation and localized structures of the three solitons, suggesting its potential extension to general $N$-soliton solutions.

\section{Conclusion}
\label{sec:con}
In this work, we develop an extension of the scaled TW-PINN framework for solving the Korteweg--de Vries (KdV) equation with different dispersion parameters using a single trained PINN model.
For one-soliton solutions, we employ the existing scaled TW-PINN framework while adapting the scaling transformation to the KdV equation.
For two-soliton solutions, the framework is extended by incorporating the propagation characteristics of individual solitons into the network formulation.
In addition, we propose soliton-informed sampling to incorporate the predicted soliton trajectories into the sampling process.

Numerical experiments show that the scaled TW-PINN accurately reproduces one-soliton propagation for different amplitudes and dispersion parameters.
For the two-soliton problem, Step 1 accurately estimates the propagation speeds of individual solitons.
To isolate the effects of the traveling-wave coordinates and SIS from speed estimation errors, we use the exact propagation speeds in Step 2.
The resulting solver accurately captures the nonlinear interactions between solitons and preserves their structures during extrapolation.
Furthermore, the three-soliton results demonstrate that TW-PINN with the same two-step training procedure captures the propagation and nonlinear interactions of three solitons.
These results suggest that the procedure may also be applicable to general $N$-soliton problems.

Future work will focus on extending the proposed framework to more general multi-soliton configurations, more general initial conditions beyond exact soliton profiles, and more complex nonlinear dispersive systems, where intrinsic solution structures can be further exploited to improve the accuracy and efficiency of physics-informed learning methods.


\providecommand{\bysame}{\leavevmode\hbox to3em{\hrulefill}\thinspace}
\providecommand{\MR}{\relax\ifhmode\unskip\space\fi MR }
\providecommand{\MRhref}[2]{%
  \href{http://www.ams.org/mathscinet-getitem?mr=#1}{#2}
}
\providecommand{\href}[2]{#2}

\end{document}